\documentclass{gtart}
\usepackage{natbib, amssymb,latexsym, amscd}
\usepackage[all]{xy}
\usepackage{graphicx}

 \newtheorem{theorem}{Theorem}[section]
 \newtheorem{cor}[theorem]{Corollary}
 \newtheorem{lemma}[theorem]{Lemma}
 \newtheorem{proposition}[theorem]{Proposition} \theoremstyle{definition}
 \newtheorem{definition}[theorem]{Definition}
 \theoremstyle{definition}
 \newtheorem{example}[theorem]{Example}
 \theoremstyle{remark}
 \newtheorem{rem}[theorem]{Remark}

\usepackage{latexsym}
\usepackage{amssymb}
\newcommand{\ben}{\begin{equation}}
\newcommand{\een}{\end{equation}}

\newcommand{\integer}{\ensuremath{{\mathbb Z}}}

\newcommand{\real}{\ensuremath{{\mathbb R}}}
\newcommand{\complex}{\ensuremath{{\mathbb C}}}
\newcommand{\quaternion}{\ensuremath{{\mathbb H}}}
\newcommand{\rational}{\ensuremath{{\mathbb Q}}}

\newcommand{\U}[1]{\ensuremath{{\mathrm U( #1 )}}}

\newcommand{\DD}{{\mathcal D}}
\newcommand{\PP}{{\mathcal P}}

\newcommand{\BB}{{\mathcal B}}

\newcommand{\EE}{{\mathcal E}}

\newcommand{\CC}{\mathcal{C}}

\newcommand{\LL}{\mathcal{L}}
\newcommand{\MM}{\mathcal{M}}

\newcommand{\OO}{\mathcal{O}}

\newcommand{\Xx}{\mathsf{X}}
\newcommand{\Gg}{\mathsf{G}}

\newcommand{\target}{\mathsf{t}}
\newcommand{\source}{\mathsf{s}}

\newcommand{\mult}{\mathsf{m}}

\newcommand{\Loop}{\mathsf{L}}

\newcommand{\To}{\longrightarrow}

\newcommand{\timests}{\: {}_{\target}  \! \times_{\source}}

\newcommand{\toparrow}[1]{\stackrel{#1}{\rightarrow}}

\newcommand{\coarse}[1] {{#1} /_{\sim}}

\newcommand{\gr}{\mathfrak}

\newcommand{\dd}{\ensuremath{\DD}}

\newcommand{\Map}{\ensuremath{{\mathrm{Map}}}}

\newcommand{\Bun}{\ensuremath{{\mathrm{Bun}}}}

\newcommand{\hol}{\ensuremath{\mathrm{hol}}}
\newcommand{\ih}{\ensuremath{\mathrm{ih}}}
\newcommand{\oh}{\ensuremath{\mathrm{oh}}}
\newcommand{\td}[1]{\tilde{#1}}
\newcommand{\BF}[1]{{\bf{#1}}}
\newcommand{\mchord}{\ensuremath{\MM\CC}}

\newcommand{\cd}{\ensuremath{\CC}}
\newcommand{\lcd}{\ensuremath{\EE\CC}}
\newcommand{\mlcd}{\ensuremath{\MM\EE\CC}}

\newcommand{\mcd}{\ensuremath{\MM\CC}}
\newcommand{\mdd}{\ensuremath{\MM\DD}}
\newcommand{\gmlcd}{\ensuremath{G \MM\EE\CC}}
\newcommand{\bargmlcd}{\ensuremath{\overline{ G \MM\EE\CC}}}
\newcommand{\gmcd}{\ensuremath{G \MM\CC}}
\newcommand{\gmdd}{\ensuremath{G \MM\DD}}

\begin{document}

\title[Orbifold String Topology]{ Orbifold String Topology}

\author[E. Lupercio, B. Uribe and M. Xicot\'encatl]{Ernesto Lupercio, Bernardo Uribe and Miguel A. Xicot\'encatl}

\address{Departamento de Matem\'{a}ticas, CINVESTAV,
     Apartado Postal 14-740
     07000 M\'{e}xico, D.F. M\'{E}XICO}
\address{Departamento de Matem\'{a}ticas, Universidad de los Andes,
Carrera 1 N. 18A - 10, Bogot\'a, COLOMBIA}
\address{Departamento de Matem\'{a}ticas, CINVESTAV,
     Apartado Postal 14-740
     07000 M\'{e}xico, D.F. M\'{E}XICO}
\email{ lupercio@math.cinvestav.mx } \email{buribe@uniandes.edu.co} 
\email{xico@math.cinvestav.mx}

\begin{abstract} In this paper we study the string topology
(\`a la Chas-Sullivan) of an orbifold. We define the string homology ring
product at the level of  the free loop space of the classifying space of an
orbifold. We study its properties and do some explicit calculations.
\end{abstract}


\maketitle

\section{Introduction}
String topology is the study of the topological properties of the free loop space $\LL M$ of a smooth manifold $M$ by the use of methods originating in quantum field string theories and in classical algebraic topology. Here $\LL M$ is by definition the space $\Map(S^1;M)$ of piecewise
smooth maps from the unit circle $S^1$ to $M$. This study was initiated by Chas and Sullivan  in their seminal paper \cite{ChasSullivan}. In this paper they defined a remarkable product $\circ$ on the homology $H_*(\LL M)$ of the loop space of a smooth manifold.

The Chas-Sullivan string product
on $H_*(\LL M)$ was only part of a very interesting
structure unveiled in their work; for example Chas and Sullivan defined a degree one map $$\Delta \colon H_*(\LL M) \to H_{*+1}(\LL M)$$ given by
$$\Delta(\sigma) = \rho_*(d\theta \otimes \sigma)$$
where $\rho\colon S^1 \times \LL M \to \LL M$ is the evaluation map and $d\theta$ is the fundamental class of $S^1$. One of the main theorems of \cite{ChasSullivan} is the following one

\begin{theorem}[Chas-Sullivan \cite{ChasSullivan}] \label{CSBValg}
The triple $$(H_*(\LL M), \circ, \Delta)$$ is a Batalin-Vilkovisky
algebra, namely
\begin{itemize}
    \item $(H_{*-d}(\LL M), \circ)$ is a graded commutative algebra.
    \item $\Delta^2=0$
    \item The bracket $$\{\alpha,\beta\} = (-1)^{|\alpha|} \Delta(\alpha \circ \beta) - (-1)^{|\alpha|} \Delta(\alpha)\circ \beta - \alpha\circ\Delta(\beta)$$ makes $H_{*-d}(M)$ into a graded Gerstenhaber algebra (namely it is a Lie bracket which is a derivation on each variable).
\end{itemize}
\end{theorem}

This establishes a striking relation between algebraic topology and recent findings in
quantum field theory and string theory \cite{BV, Getzler}.

Cohen and Jones \cite{CohenJones} discovered that a
very rich part of this structure was available at a more homotopy-theoretic
level and reinterpreted the BV-algebra structure in terms of an action of the
cactus operad on a certain prospectrum associated to $M$. They showed moreover
that the Chas-Sullivan string product was the natural product in the Hochschild
cohomology interpretation of the homology of the loop space of $M$ \cite{Jones}. Cohen and Godin \cite{CohenGodin} studied interactions with the study of the homology of moduli spaces of Riemann
surfaces, establishing a direct connection to topological quantum field theories. Cohen and Godin used the concept of Sullivan chord diagram in their work. Cohen, Jones and Yan \cite{CohenJonesYan} provided more explicit calculations of the product by the careful use of the spectral sequence associated to the  fibration $$ \Omega M \longrightarrow \LL M \longrightarrow M$$ induced by the evaluation map. In particular they computed the Chas-Sullivan product on the homology of the free loop space of spheres and complex projective spaces.

In this paper we generalize several of the fundamental results of string topology by showing that they remain true if we replace the manifold $M$ by an orientable orbifold $\Xx=[M/G]$, where $G$ is a finite group acting by orientation preserving diffeomorphisms on $M$. More precisely the following theorem is  the main result in this paper and can be seen as a generalization of Theorem \ref{CSBValg} to the orbifold context.

\begin{theorem}\label{OrbSTalg} Let $\Xx=[M/G]$ be an orientable orbifold, then $$A_{\Loop \Xx}:=H_*(\LL(M\times_G EG); \rational)$$ has the structure of a Batalin-Vilkovisky algebra.
\end{theorem}

This BV-algebra can be identified in two extreme cases:
\begin{itemize}
    \item When $G=\{1\}$ and for arbitrary $M$ then $A_{\Loop \Xx}$ coincides with the Chas-Sullivan BV-algebra.
    \item When $M=\{m_0\}$ is a single point and for arbitrary finite $G$ then $A_{\Loop \Xx}$ is isomorphic to the center of the group algebra of $G$.
\end{itemize}

This paper is organized as follows. In section \ref{section free loop space} we define a topological groupoid that we call \emph{the loop groupoid} $\Loop\Xx$ \cite{LupercioUribeLoopGroupoid} of an orbifold groupoid $\Xx$ and prove that is satisfies the basic property $$B\Loop\Xx \simeq \LL B \Xx$$ where $B$ if the geometric realization functor from groupoids to topological spaces \cite{Segal1}.  In section \ref{The orbifold string  product.} we define the product in $H_*(\LL(M\times_G EG); \rational)$ that eventually  will be part of the BV-algebra structure. In section \ref{Discrete torsion} we show that there is a version of the ring of section \ref{The orbifold string  product.} twisted by discrete torsion which is analogous to the twistings in equivariant $K$-theory of Adem and Ruan \cite{AdemRuan}.  In section \ref{Operadic structure} we prove Theorem \ref{OrbSTalg} which is our main result. For this purpose we introduce and operad of Sullivan chord diagrams which we call the \emph{chord diagram operad} and show that it is equivalent to the cactus operad.
 Finally in Section \ref{Computations} we perform some explicit computations, for example we show how to use the methods of this paper to compute the Chas-Sullivan product on Lens spaces.

We would like to thank Alejandro Adem, Carl-Friedrich B\"odigheimer, Ralph Cohen, Dan Freed, Sam Gitler, Maxim Kontsevich, Jacob Mostovoy, Mainak Poddar, Antonio Ramirez, Yongbin Ruan, Graeme Segal, Jim Stasheff and Dennis Sullivan for relevant conversations regarding this work.

We acknowledge the support of CONACYT (the first and third authors) and of the Universidad de los Andes (the second author).

\section{The free loop space of the classifying space of an orbifold} \label{section free loop space}

\subsection{Orbifolds.}

For an introduction to the theory of orbifolds, we refer the reader to \cite{Moerdijk2002} or \cite{AdemLeidaRuan}. Here we just recall the basic definitions we will use in this paper.

Following Moerdijk \cite{Moerdijk2002}, we will use groupoids to study orbifolds. A groupoid is a (small) category in which each arrow is an isomorphism. Whenever we have a groupoid $\Gg$ we will denote by $\Gg_0$ the set of all its objects and by $\Gg_1$ the set of all its arrows.  We will denote by $\Gg_1 \timests \Gg_1$ the subset of $\Gg_1 \times \Gg_1$ consisting pairs of arrows so that the target of the first equals the source of the second. We will write $1_x$ to denote the identity arrow for the object $x$ in $\Gg_0$.

We will denote the structure maps by:
      $$\xymatrix{
        \Gg_1 \timests \Gg_1 \ar[r]^{\ \ \ m} & \Gg_1 \ar[r]^i &
        \Gg_1 \ar@<.5ex>[r]^s \ar@<-.5ex>[r]_t & \Gg_0 \ar[r]^e & \Gg_1
      }$$
where $s$ and $t$ are the source and the target maps,  $m$  is the multiplication map given by $$m(\alpha, \beta) = \beta \circ \alpha,$$ $e$ is the unit map defined by $$e(x) = 1_x,$$ and $i$ is the inverse map $$i(\alpha) = \alpha^{-1}.$$

\begin{definition} A \emph{Lie groupoid} is a groupoid $\Gg$ for which, in addition, $\Gg_0$ and $\Gg_1$ are Hausdorff smooth manifolds, and the structure maps $s$, $t$, $m$, $e$ and $i$ are all smooth. We will also require   $s$ and $t$ to be submersions.

We define several kinds of Lie groupoids.

\begin{itemize}
\item A Lie groupoid $\Gg$ is said to be a \emph{proper} groupoid if the map $(s,t) \colon \Gg_1 \longrightarrow \Gg_0\times \Gg_0$ is proper.
\item A Lie groupoid $\Gg$ is said to be a \emph{foliation} groupoid if for every object $x \in \Gg_0$ its stabilizer $\Gg_x:=\{\alpha \in \Gg_1 | s(\alpha) = t(\alpha) = x\}$ is finite.
\item A proper foliation  Lie groupoid $\Gg$ will be called an \emph{orbifold} groupoid.
\end{itemize}
\end{definition}

\begin{example} \label{translation} Let $G$ be a finite group acting smoothly from the right on a smooth manifold $M$. Then the \emph{translation groupoid} $M \rtimes G$ is defined by $(M \rtimes G)_0 = M$, $(M \rtimes G)_1=M \times G$, $s(x,g)=x$, $t(x,g)=xg$, $m((x,g),(xg,h))=(x,gh)$, $i(x,g)=(x,g^{-1})$ and $e(x)=(x,1_G)$.
\end{example}

Moerdijk shows \cite{Moerdijk2002} that the orbifold groupoids together with smooth functions form a category, and calls a morphism in this category a \emph{Morita equivalence} if it is fully faithful and essentially surjective in a an appropriate differential-geometric sense. He then defines the category of orbifolds as the localization of the category of orbifold groupoids at the class of Morita equivalences. If $M$ and $G$ are as in Example \ref{translation}, we will write $[M/G]$ for the translation groupoid $ M \rtimes G$ viewed as an object of the orbifold category. Orbifolds of this form are called \emph{`global quotient orbifolds'}. For simplicity, we will restrict our attention to such global quotients orbifolds throughout the paper. A reader not familiar with orbifolds may just think of the data of such global quotient orbifold as equivariant data, keeping in mind that there is a notion of equivalence which might identify the data of an action of a group $G$ on $M$ with those of an action of a different group $\Gamma$ on a different manifold $N$.

Segal in \cite{Segal1} defines a functor from the category of all small categories to the category of topological spaces  $$ B:\mathbf{Cat} \longrightarrow \mathbf{Top}$$ that he calls the classifying space functor. This functor sends transformations of functors to homotopies of continuous mappings.

The \emph{classifying space of the orbifold $\Xx=[M/G]$} is defined as the classifying space $B(M \rtimes G)$ of the category $M \rtimes G$. It is known \cite{Moerdijk2002} that $B(M \rtimes G)$  depends only on the orbifold and not on the particular group action used to represent it. For this reason we will denote it also by $B\Xx$.

Let $\Delta^n$ be the standard $n$-simplex. If we denote by $\Xx_n$ the space of ordered $n$-tuples of arrows in $M \rtimes G$ such that the target of each arrow is equal to the source of the following one, then  it is not hard to see that $$\Xx_n \cong M \times G \times \cdots \times G$$ where $G$ appears exactly $n$ times. The classifying space is defined as a quotient
\begin{equation} \label{clasificante}
 B\Xx := \bigsqcup_n (\Xx_n  \times \Delta^n) / \sim
\end{equation}
where $\sim$ is an equivalence relation defined in \cite{Segal1}.

The space $B\Xx$ is homotopy equivalent to  the Borel construction of the $G$-manifold $M$, namely
\begin{equation}\label{Borel}
B \Xx := B(M\rtimes G) \simeq M_G := M \times_G EG,
\end{equation}
where $EG$ is a contractible space where $G$ acts freely and $M \times_G EG := (M \times EG)/G$ is the quotient by the diagonal action. For a proof of this fact we refer the reader to \cite{Segal1}. We will write $[m,\xi]$ to denote the equivalence class of $(m,\xi) \in M  \times EG$ under this diagonal action. We will also write $[Y,T]$ to denote the projection of the set $Y \times T \subseteq M \times EG$ to $B\Xx$

Let $Y$ and $Z$ be two $G$-spaces, and let $$e:Z\to Y$$ be an equivariant map between them (i.e. $e(zg)=e(z)g$ for all $z$ and $g$). The assignment that associates to a $G$-space its translation groupoid is functorial. In other words $e$ induces a morphism of groupoids which we call again $e$. We will denote by $$|e|:Z\times_G EG \longrightarrow Y\times_G EG$$ the corresponding map of classifying spaces.

\subsection{The loop space as a classifying space for the loop orbifold.} \label{loop orbifold}

Denote by $\LL B \Xx = \LL (M_G)$ the free loop space of $B \Xx$, then $\LL B \Xx$ depends
only on $\Xx$. In this section we will find a topological groupoid $\Loop \Xx$ with the property that $B \Loop \Xx \simeq \LL B \Xx$.

We find the following notation useful. We set $$\PP^x_y(M) \colon = \{ \gamma
\colon [0,1] \to M | \gamma(0)= x\ \mbox{and}\  \gamma(1)=y\}$$ that is to say
the space of all piecewise smooth paths going from $x \in M$ to $y \in M$.

We define the \emph{loop groupoid}
$\Loop \Xx$  of a global
orbifold $\Xx$ (for a definition in the general case see \cite{LupercioUribeLoopGroupoid} cf. \cite{BridsonHaefliger}).  For this we consider the space
$$\PP_G(M):= \bigsqcup_{g \in G} \PP_g(M)\times \{g\}$$ where
$$\PP_g(M) := \{ \gamma \colon [0,1] \to M | \gamma(0) g = \gamma(1) \} = \bigsqcup_{x \in M} \PP^x_{xg}(M)$$
together with the $G$ action given by
\begin{eqnarray*}
\PP_G(M) \times G & \to & \PP_G(M) \\
((\gamma,g),h) & \mapsto & (\gamma_h , h^{-1}gh)
\end{eqnarray*}
where $\gamma_h(m) : = \gamma(m)h$.
\begin{eqnarray} \includegraphics[height=2.0in]{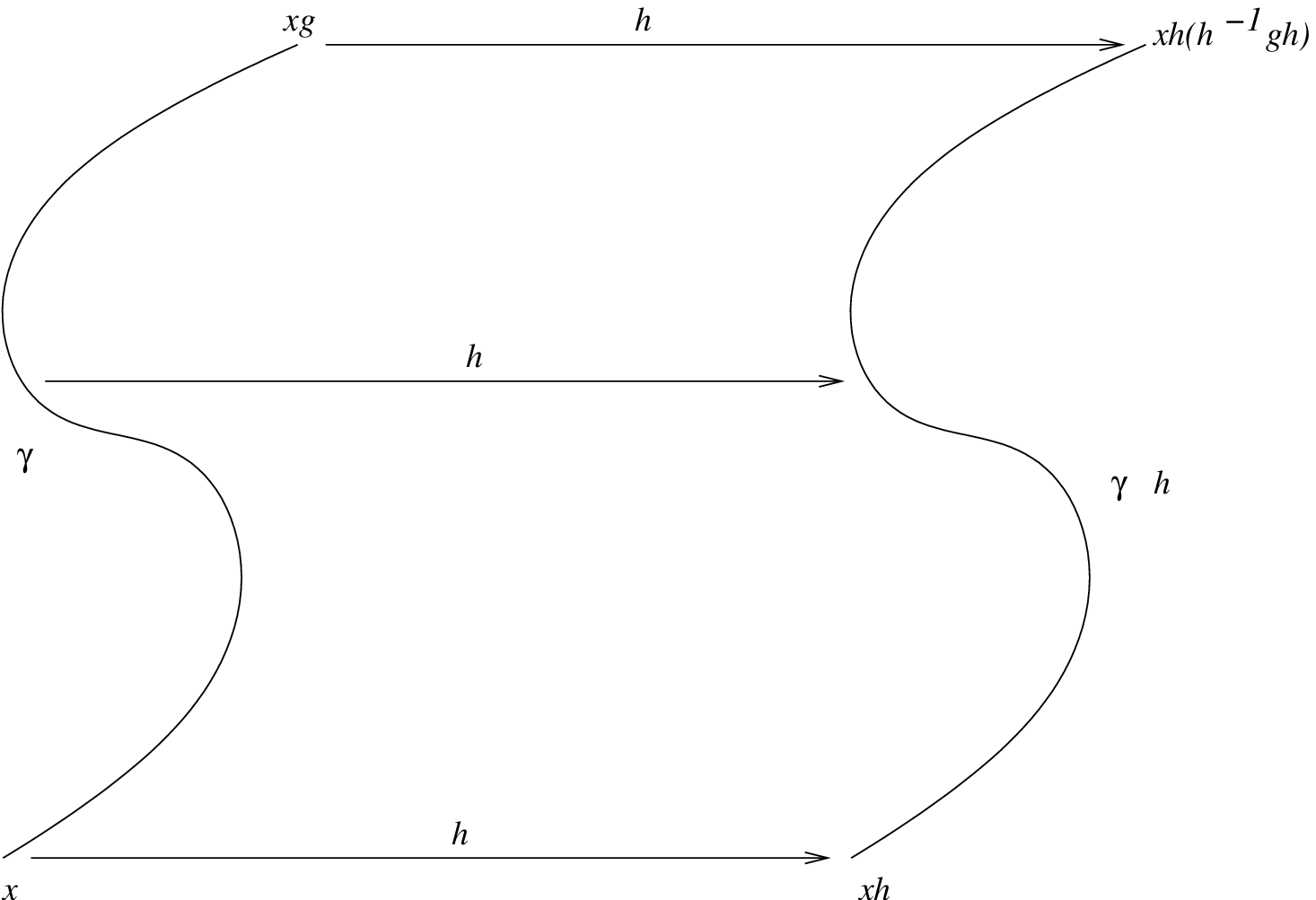} \label{chord}
\end{eqnarray}
Then  we define the
loop groupoid as $$\Loop \Xx : = \PP_G(M) \rtimes G.$$

The following result describes the relation between $\LL  B \Xx$ and $\Loop \Xx$

\begin{theorem} \label{theorem loop orbifold}
There is a canonical map
$$ \tau \colon \LL B \Xx\To B \Loop \Xx $$
that induces a weak homotopy equivalence.
\end{theorem}

\begin{proof}
We will construct two Serre fibrations over $B \Xx$.
\begin{itemize}
\item  Consider $B \Loop \Xx$.

Define a morphism of groupoids
$$\widetilde{ev}_0 \colon \Loop \Xx \to \Xx$$
induced by the equivariant map of $G$-spaces $$ev_0 \colon  \PP_G(M) \longrightarrow M$$ given by evaluation at $0$, $$ev_0(\gamma,g) : = \gamma(0).$$

 This morphism induces a map at the level of classifying spaces
$$| ev_0 | \colon B \Loop \Xx \to B \Xx.$$
If we interpret the classifying spaces in terms of the Borel construction as in formula \ref{Borel} we have $B
\Loop \Xx = \PP_G(M) \times_G EG$ and $B \Xx = M \times_G EG$. For a point $z
\in B \Xx$ with $z = [m, {\xi}]$, the following holds
$$|ev|^{-1}(z) = \left[ \PP_G^m(M) \times \{\xi\}\right]$$
where $$\PP_G^m(M) : = \bigsqcup_{g \in G} \PP_g^m(M) \times \{ g\}$$ with
$$\PP_g^m(M) = \{ \gamma \in \PP_g(M) | \gamma(0) =m \}.$$


\item On $\LL B \Xx$.

Take the map $$\epsilon_0 \colon \LL B \Xx \To B \Xx$$ which evaluates a free
loop at $0$, i.e. for $\sigma \colon S^1 \to B \Xx$ then $\epsilon_0(\sigma) : =
\sigma(0)$. Then $$\epsilon_0^{-1}(z) = \Omega_z(B\Xx) : = \PP^z_z(B\Xx)$$ is the space of loops
based at $z$.
\end{itemize}

Now let's define the map $\tau$. Consider the fixed $z=[m,\xi]$ as above and for $\sigma \in \LL B \Xx$, lift it to $\tilde{\sigma}$ making the following
diagram commutative

$$
  \xymatrix{
    [0,1] \ar[r]^(.4){\tilde{\sigma}} \ar[d]_{\exp(2\pi\_ \ )} & M \times EG \ar[d]^p \\
    S^1 \ar[r]^(.4)\sigma & M \times_G EG
  }
$$
such that $\tilde{\sigma}(0)=(m, {\xi})$ (the construction follows from
the fact that the map $p$ is a $G$-principal bundle and $G$ is finite). Since $G$ acts freely on $EG$ there exists a unique element $k$ in $G$ such that $\tilde{\sigma}(0)k = \tilde{\sigma}(1)$. Define $\tau$ in the following way
$$\tau(\sigma) \colon = \left[  (\pi_1 \circ \tilde{\sigma}, k ), {\xi}  \right ] \in B \Loop \Xx,$$
where $\pi_1:M \times EG \to M$ is the projection on the first coordinate.
From the definition of $\tau$ it follows that it is well defined and that $\pi_1 \circ \tilde{\sigma} \in \PP_k(M)$. Moreover the following diagram is commutative
$$
\xymatrix{
\LL B \Xx \ar[dr]_{\epsilon_0}  \ar[rr]^\tau & & B \Loop \Xx \ar[dl]^{|ev_0|}\\
 & B \Xx. &
}
$$

Let's denote by $\tau_z : = \tau|_{\epsilon_o^{-1}(z)}$, then

\begin{lemma}
The map
$$\tau_z \colon \epsilon_0^{-1}(z) \To |ev_0|^{-1}(z) $$
is a homotopy equivalence.
\end{lemma}

\begin{proof} From the definition of $\tau$ it is clear that $\tau_z$ is surjective. Let's now check the homotopy type of the inverse image of a point.
  Recall from above that the map
$\tau_z$ goes from $\Omega_z(M\times_G EG)$ to $\left[ \PP_G^m(M) \times
\{\xi\}\right]$. Take $(\gamma,g) \in  \PP_G^m(M)$.  From the definition of
$\tau$ above it follows that
$$\tau_z^{-1}\left( \left[ \left( (\gamma,g), {\xi} \right) \right]\right) \cong \PP^{{\xi}}_{g^{-1}{\xi}}(EG)$$
where $\PP^{{\xi}}_{g^{-1}{\xi}}(EG)$ stands for the paths in
$EG$ that go from ${\xi}$ to $(g^{-1}{\xi})$.

The space $ \PP^{{\xi}}_{g^{-1}{\xi}}(EG)$ is independent of the choice of
representative in $\left[ \left( (\gamma,g), {\xi} \right) \right]$.

As the space $\PP^{{\xi}}_{g^{-1}{\xi}}(EG)$ is contractible
then it follows that $\tau_z$ induces a homotopy equivalence.
\end{proof}

As $\tau$ induces a homotopy equivalence on the fibers of the Serre fibrations
given by $\epsilon_0$ and $|ev_0|$, then the theorem \ref{theorem loop orbifold}
follows from a theorem of Dold \cite{Dold}. Hence $\tau$ induces a weak homotopy
equivalence between $\LL B \Xx$ and $B \Loop \Xx$.
\end{proof}

\subsection{Circle action} \label{circle action}

We have seen that the map $ \tau \colon \LL B \Xx\To B \left(\Loop
\Xx \right)$ is a weak homotopy equivalence, and it is natural to
wonder whether the equivalence is $S^1$-equivariant. The answer
turns out to be negative as we will see shortly.

There is a natural action of $S^1$ onto $\LL B \Xx$ by rotating the
loop, but the action does not get carried into $B \Loop \Xx$ via
$\tau$. The reason is the following, the loop orbifold $\Loop \Xx$ comes
provided with a natural action of the orbifold $[\real /
\integer]$ which is a \emph{stack} model for the circle. The action of $\real $ into the
orbifold loops of $\PP_G(M)$ is the obvious one, the map gets
shifted by the parameter in $\real$. The subtlety arises here,
once we act on the orbifold loop by $1\in \real$, we do not end up
with the orbifold loop from the beginning, but instead we get one
that is related to the initial one via an arrow of the loop
orbifold category. This arrow in the loop orbifold is where $1 \in
\integer$ gets mapped. By the way, precisely this fact was the one that
allowed us to define the loop orbifold in a non trivial way,
namely a loop on the orbifold was not a map from the circle to the
orbifold, but a functor from $[\real / \integer]$ to the orbifold.

More accurately, to define the action of $[\real
/ \integer] $ on $\Loop \Xx = [\PP_G(M) / G]$ we first define an
action of $\real$ on $\PP_G(M)$ in the natural way, namely, take
$\gamma \in \PP_k(M)$ and $s\in \real$ and define $$(s \cdot \gamma) (t):=
\gamma_s(t) = \gamma(t+s -\lfloor t+s \rfloor)k^{\lfloor t+s \rfloor}$$
where $\lfloor \cdot \rfloor$ is the least integer function. Then for each $(\gamma, k) \in
\PP_G(M)$ and $1 \in \integer$ we choose the arrow of $\Loop \Xx$
that relates the orbifold loops $(\gamma, k)$ and $(\gamma_1,k)$,
this is  the arrow $((\gamma,k),k) \in \PP_G(M) \times G$. The
source of $((\gamma,k),k)$ is $(\gamma,k)$ and the target is
$(\gamma\cdot k,k) = (\gamma_1,k)$ the loop shifted by $1$.

Using the construction of section \ref{loop orbifold} we have that
$$\tau(\sigma) \colon = \left[ \left( (\pi_1 \circ \tilde{\sigma}, k ),
{\xi} \right) \right ],$$ and denote $\pi_1 \circ
\tilde{\sigma} $ by $ \gamma$. For $s \in \real $,
$$ \tau (s \cdot \sigma) =
\left[ \left( \gamma_s, k ), {\xi} \right) \right ]$$ and
$1 \cdot \sigma = \sigma$, but $\tau (1 \cdot \sigma) \neq
\tau(\sigma)$. Instead $\tau (1 \cdot \sigma)$ and $ \tau(\sigma)$
are related by an arrow.

Nevertheless, if we take the coarse moduli space of $\Loop \Xx$ (that we will write $\coarse{\Loop \Xx} = \PP_G(M)/G$),
the map induced by $\tau$ is $S^1$-equivariant. For in
$\coarse{\Loop \Xx} = \PP_G(M)/G$  the elements
$\tau (1 \cdot \sigma)$ and $ \tau(\sigma)$ become by definition the same. Then
we can conclude

\begin{lemma} The space $\coarse{\Loop \Xx} = \PP_G(M)/G$ has a natural $S^1$
action and the map  $$\tilde{\tau} : \LL B \Xx \To \coarse{\Loop
\Xx} = \PP_G(M)/G$$ which is the composition of $\tau$ with the
projection $B \Loop \Xx \to \coarse{\Loop \Xx}$,
is $S^1$-equivariant.
\end{lemma}

\begin{cor} \label{corollary circle action}
The map $\tilde{\tau}$ induces an isomorphism in homology
$$\tilde{\tau}_* : H_*(\LL B \Xx; \rational) \stackrel{\cong}{\To}
H_*(\coarse{\Loop \Xx};\rational),$$ and in equivariant homology
$$\tilde{\tau}_* : H^{S^1}_*(\LL B \Xx; \rational) \stackrel{\cong}{\To}
H^{S^1}_*(\coarse{\Loop \Xx};\rational)$$
\end{cor}
\begin{proof}
As $\tau$ is a weak homotopy equivalence, then
$$\tau _* : H_*(\LL B \Xx ; \integer) \stackrel{\cong}{\To} H_*(B \Loop \Xx
; \integer),$$ and as the group $G$ is finite then
$$ \tilde{\tau}_* : H_*(B \Loop \Xx ; \rational) \stackrel{\cong}{\To}
H_*(\coarse{\Loop \Xx}; \rational).$$ The second isomorphism follows from
the isomorphism of spectral sequences with real coefficients
associated to the each of the following fibrations
$$
\xymatrix{ \LL B \Xx \times_{S^1} ES^1 \ar[rr]\ar[dr] & &
\coarse{\Loop \Xx} \times_{S^1} ES^1 \ar[dl]\\
& BS^1. &}
$$

\end{proof}

\subsection{Cyclic equivariant loops}

There is an alternative description of $\PP_g(M)$ that although essentially
obvious nevertheless relates it to some models that have been studied before.

Given an element $g\in G$ it generates a cyclic group $<g> \subseteq G$.  Let
$m$ be the order of $g$ in $G$. Then there is a natural injective morphism of
groups $$ \zeta \colon <g> \to S^1$$ given by $\zeta(g) = \exp( 2 \pi i / m)$.

We define the space $\LL_g M$ of \emph{$g$-equivariant loops} in $M$ to be the
subspace of $\LL M := \Map(S^1 ; M)$ of loops $\phi$ satisfying the following
equation for every $z \in S^1$: $$ \phi(z \cdot \zeta(g))=\phi(z)\cdot g.$$

The space of \emph{cyclic equivariant loops} of $M$ is defined to be  simply
$$\LL_G M := \bigsqcup_{g \in G} \LL_g M \times\{g\}.$$ It is, again, naturally
endowed with a $G$-action $((\phi,h);g) \mapsto (\phi g, g^{-1}hg)$.

The natural restriction map $$\Psi\colon \LL_g(M) \To \PP_g(M)$$ given by $$
\gamma(t) = \phi(\exp(2 \pi i t / m) )=  \phi(\zeta(g)^t)$$ is a diffeomorphism,
and moreover it induces a $G$-equivariant diffeomorphism $$\Psi\colon \LL_G(M)
\To \PP_G(M).$$

We conclude this subsection by pointing out that as a consequence of these
remarks we have the following equality $$ \LL_G(M) \times_G EG \simeq \LL(M_G) =
\LL B \Xx.$$

\subsection{Principal bundles.}\label{principal bundles}

In this paragraph we consider $G$-principal bundles on $S^1$ and their relation
to the various models of the loop orbifold. We are interested in the category
of $G$-principal bundles $\pi\colon Q \to S^1$ over $S^1$ endowed with a marked
point $q_0 \in Q$ so that $\pi(q_0)=0 \in S^1$, and such that $\pi$ is a local isometry.

Whenever we have such a pair $(Q,q_0)$ we have a well-defined  lift $\tilde{e}
\colon [0,1] \to Q$, $\tilde{e}(0)=q_0$, of the exponential map $e \colon [0,1] \to S^1$ given by $t
\mapsto \exp(2 \pi i t)$, making the following diagram commutative:
$$
\xymatrix{
[0,1] \ar[dr]_{e}  \ar[rr]^{\tilde{e}} & & Q \ar[dl]^{\pi}\\
 & S^1. &
}
$$
Since $\tilde{e}(0)$ and $\tilde{e}(1)$ belong to $\pi^{-1}(0)$ there is a $g \in
G$ so that $$\tilde{e}(1) = \tilde{e}(0) \cdot g.$$ We well call this $g\in G$
the \emph{holonomy} of $Q$.

The isomorphism classes of $G$-principal bundles with a marked point are
classified by their holonomy, for the set $\Bun_G(S^1)$ of such classes is given
by $$\Bun_G(S^1)= \pi_1 BG = G.$$ The following proposition is very
easy.

\begin{proposition} The natural action of $G$ on $\Bun_G(S^1)$ under the
holonomy isomorphism $\hol \colon \Bun_G(S^1) \To G$ becomes the action of $G$
on $G$ by conjugation.
\end{proposition}

This proposition can be slightly generalized as follows. Consider now the space
$\Bun_G(S^1,M)$ of isomorphism classes of $G$-equivariant maps from a principal
$G$-bundle $Q$ over the circle to $M$. This space has a natural $G$-action defined as follows.
If $Q_g$ denotes the principal bundle with holonomy $g$ then the pair
$$[(\beta\colon Q_g \to M); k] \in \Bun_G(S^1,M) \times G$$ gets mapped by
conjugation to $$(\beta_k \colon Q_{k^{-1}gk} \to M) \in \Bun_G(S^1,M).$$

\begin{proposition} The loop orbifold $\Loop \Xx = [\PP_G(M)/ G]$ is isomorphic
to the orbifold $[\Bun_G(S^1,M) / G]$, and therefore $$\Bun_G(S^1,M)\times_G EG
\simeq \LL(M_G).$$
\end{proposition}
\begin{proof} It is enough to give a $G$-equivariant diffeomorphism
$$\Bun_G(S^1,M)\To \PP_G(M),$$ this can be achieved by the following formula $$
(\beta \colon Q_g \to M ) \mapsto \gamma = \beta \circ \tilde{e}. $$ Since
$\tilde{e}(1) = \tilde{e}(0) \cdot g$, then $\gamma(1) = \gamma(0) \cdot g.$
\end{proof}

To finish this section let us define $\Bun_g(S^1,M)$ to be the space of
isomorphism classes of $G$-equivariant maps from a principal $G$-bundle $Q_g$
with holonomy $g$ to $M$. Then we have that $$\Bun_G(S^1,M)= \bigsqcup_{g\in G}
\Bun_g(S^1,M),$$ and in fact $$\Bun_g(S^1,M) \cong \PP_g(M).$$

\section{The orbifold string product}\label{The orbifold string product.}

\subsection{The definition of the orbifold string ring.} We will suppose at the outset of what follows that $M$ is oriented and that $G$ acts in an orientation preserving fashion. In this section we define the string product in the homology of $\LL B
\Xx$. Let us start by defining a composition-of-paths map: $$ \star \colon
\PP_g(M) \: {}_{\epsilon_1} \! \! \times_{\epsilon_0} \PP_h(M) \To \PP_{gh}(M)$$ Here the
map $\epsilon_t \colon \PP_k(M) \to M$ is the evaluation map at $t$ given by
$\gamma \mapsto \gamma(t)$. Therefore $\PP_g(M) \: {}_{\epsilon_1} \! \! \times_{\epsilon_0} \PP_h(M)$ is the space of pairs of paths
$(\gamma_0, \gamma_1)$ so that
$\gamma_0(1)=\gamma_1(0)=:\epsilon_\infty(\gamma_0,\gamma_1)$, that is to say the end
of $\gamma_0$ is the beginning of $\gamma_1$. The map $\star$ is given by the
obvious formulas of concatenation:
\begin{eqnarray*}
(\gamma_0 \star \gamma_1 )(t) := \left\{
\begin{array}{ccc}
\gamma_0 (2t) & \mbox{for} &  0 \leq t \leq 1/2 \\
\gamma_1 ( 2t-1) & \mbox{for} & 1/2 < t \leq 1
\end{array} \right.
\end{eqnarray*}
Notice that the following diagram is a cartesian square:
$$
\xymatrix{
 \PP_g(M) \: {}_{\epsilon_1} \! \! \times_{\epsilon_0} \PP_h(M) \ar[r]^(.52)j \ar[d]_{\epsilon_\infty} & \PP_g(M)\times
 \PP_h(M)\ar[d]^{\epsilon_1 \times \epsilon_0} \\
 M \ar[r]^\Delta & M\times M
 }
 $$
where $j$ is the natural inclusion and $\Delta$ is the diagonal map. Following
Cohen-Jones \cite[Section 1]{CohenJones} we observe that such a pullback square allows one to construct a
Pontrjagin-Thom map $$\tilde{\tau} \colon \PP_g(M) \times  \PP_h(M) \To
(\PP_g(M) \: {}_{\epsilon_1} \! \! \times_{\epsilon_0} \PP_h(M))^{TM}$$ where $(\PP_g(M) \: {}_{\epsilon_1} \! \! \times_{\epsilon_0} \PP_h(M))^{TM}$ denotes the Thom space of the
pullback bundle $\epsilon_\infty^*(TM)$, which is the normal bundle of the
embedding $j$. Here we remind the reader that the normal bundle of the embedding
$\Delta$ is $TM$.

Let us denote by $(\PP_{gh}(M))^{TM}$ the Thom space of the bundle
$\epsilon_{1/2}^*(TM)$ where $\epsilon_{1/2} \colon \PP_{gh}(M) \to M$. The map
$\star$ induces a map of Thom spaces $$\tilde{\star} \colon  (\PP_g(M) \: {}_{\epsilon_1} \! \! \times_{\epsilon_0} \PP_h(M))^{TM} \To (\PP_{gh}(M))^{TM}.$$ It
is immediate to verify that the following diagram is commutative
\begin{eqnarray} \label{diagram of string product}
\xymatrix{
   \PP_g(M) \times \PP_h(M) \ar[r]^(.45){\tilde{\tau}} \ar[d]_{\epsilon_1 \times
   \epsilon_0} & (\PP_g(M) \: {}_{\epsilon_1} \! \! \times_{\epsilon_0} \PP_h(M))^{TM}
   \ar[d]^{\epsilon_\infty} \ar[r]^(.6){\tilde{\star}} & (\PP_{gh}(M))^{TM} \ar[d]^{\epsilon_{1/2}} \\
   M\times M \ar[r]^{\tau} & M^{TM} \ar[r]^{=} & M^{TM}.
   }
\end{eqnarray}
From this we can consider the composition
$$
\circ \colon H_p(\PP_g(M)) \otimes H_q(\PP_h(M)) \toparrow{\times}
H_{p+q}(\PP_g(M)\times\PP_h(M)) \toparrow{\tilde{\tau}_*} $$ $$H_{p+q}((\PP_g(M) \: {}_{\epsilon_1} \! \! \times_{\epsilon_0} \PP_h(M))^{TM})\toparrow{\tilde{\star}_*}
H_{p+q}((\PP_{gh}(M)^{TM}) \toparrow{\tilde{u}_*} H_{p+q-d}(\PP_{gh}(M)), $$
where $\tilde{u}_*$ is the Thom isomorphism. By considering the direct sum over
all elements $g\in G$
we obtain the map
\begin{eqnarray}
 \circ \colon H_p(\PP_G(M)) \otimes H_q(\PP_G(M)) \To H_{p+q-d}(\PP_G(M)). \label{circ product}
\end{eqnarray}

We will call $\circ$ the \emph{$G$-string product}, and $H_*(\PP_G(M))$ the
\emph{$G$-string ring} of the $G$-manifold $M$.

This ring depends on the $G$-manifold $M$ and not only on the orbifold
$\Xx=[M/G]$. To obtain a \emph{bona fide} orbifold invariant we proceed in the
following manner.

We will define maps:
\begin{itemize}
   \item[i)] The transfer map $\theta_* \colon H_*(\LL B \Xx) = H_*(\PP_G(M)\times_G EG) \to H_*(\PP_G(M))$ which is the composition in homology of the transfer of the finite covering $\PP_G (M)\times EG \to \PP_G(M)\times_G EG$  with the projection $\PP_G(M) \times EG \to \PP_G(M)$ (see \cite[Chapter 4]{Adams}).
   \item[ii)] The projection map $\sigma_* \colon H_*(\PP_G(M)) \to H_*(\PP_G(M) \times_G EG) = H_*(\LL B \Xx)$ which is the composition in homology of the maps $$\PP_G(M) \to \PP_G(M) \times EG \to \PP_G(M)\times_G EG.$$
\end{itemize}
and then we will define the \emph{orbifold string product} $\circ$ as follows $$
\circ \colon H_p(\LL B \Xx) \otimes H_q(\LL B \Xx) \stackrel{\theta_* \otimes
\theta_*}{\longrightarrow} H_p(\PP_G(M)) \otimes H_q(\PP_G(M)) $$ $$\toparrow{\circ}
H_{p+q-d}(\PP_G(M)) \toparrow{\sigma_*} H_{p+q-d}(\LL B \Xx).$$

We will call the ring thus obtained $$\circ \colon H_*(\LL B \Xx) \otimes
H_*(\LL B \Xx) \To H_*(\LL B \Xx)$$ the \emph{ orbifold string ring} of $\Xx$.
Notice that the degree of $\circ$ is $-d$ whenever the dimension of $\Xx$ is
$d$.

Our next claim is that the string homology of the orbifold, it is
indeed an invariant on the orbifold and not of the particular choice of groupoid that we made.

Namely, we will now prove that for a functor $F:[N/\Gamma] \to [M/G]$
between two orbifolds $[N/\Gamma]$ and $[M/G]$ that defines a Morita equivalence (for
a definition of Morita equivalence see \cite{Moerdijk2002}), then
their string homologies are isomorphic. Using the definition it is not hard to check that $F=(f, \rho)$ is a pair consisting of a surjective
homomorphism $\rho : \Gamma \to G$ and a $\rho$-equivariant,
surjective, and local diffeomorphism $f: N \to M$ such that the following
diagram is a cartesian square:
$$\xymatrix{
N \times \Gamma \ar[rr]^{f\times \rho} \ar[d]^{\source \times
\target} & & M
\times G \ar[d]^{\source \times \target}\\
N\times N \ar[rr]_{f \times f} &  & M\times M }.$$

\begin{proposition}
Let $F :[N/\Gamma]\to [M/G]$ be a Morita equivalence. Then the
induced homomorphism in string homologies
$$H_*(\LL N_\Gamma; \rational ) \To H_*(\LL M_G ; \rational)$$
is an isomorphism.
\end{proposition}
\begin{proof}
In \cite{LupercioUribeLoopGroupoid} it is proven that the induced
map on loop orbifolds
$$\Loop F : \Loop [N/\Gamma] \to \Loop[M/G]$$
is a Morita equivalence. Then, their classifying spaces are
homotopically equivalent and therefore the induced homomorphism in
homologies is an isomorphism. We are left to prove the fact that
the string homology structures are the same.

As $\Loop [N/\Gamma] = [\PP_\Gamma (N) / \Gamma]$ and $\Loop [M/G]
= [\PP_G (M) / G]$ the induced map on the coarse moduli spaces
$$ \tilde{F} : \PP_\Gamma (N) / \Gamma \stackrel{\cong}{\to} \PP_G (M) / G$$
is a homeomorphism. Following the definition of the circle action
given in section \ref{circle action}, one can see that the map
$\tilde{F}$ is $S^1$-equivariant. Therefore the induced degree
shifting action in homologies $H_n \to H_{n+1}$ is the same. Now
let's see that the string product agrees.

For $\alpha, \beta \in \Gamma$ and $a = \rho(\alpha), b =
\rho(\beta)$ consider the following diagram as in (\ref{diagram of
string product}) induced by $F$.
$$
\xymatrix{
   \PP_\alpha(N) \times \PP_\beta(N) \ar[r] \ar[d] & (\PP_\alpha(N)
   \: {}_{\epsilon_1}  \!  \! \times_{\epsilon_0} \PP_\beta(N))^{TN}
   \ar[d] \ar[r] & (\PP_{\alpha \beta}(N))^{TN} \ar[d] \\
\PP_a(M) \times \PP_b(M) \ar[r]& (\PP_a(M) \: {}_{\epsilon_1}
\! \! \times_{\epsilon_0} \PP_b(M))^{TM}
    \ar[r]
     & (\PP_{ab}(M))^{TM}
   }$$

The diagram is trivially equivariant and as $TN \cong f^*TM$  then the diagram is commutative. Therefore  once
one take invariants, it produces an isomorphism of string products

$$\xymatrix{ H_p(\PP_\Gamma(N))^\Gamma \otimes H_q(\PP_\Gamma(N))^\Gamma
\ar[rr]^(.58)\circ \ar[d]_\cong   & & H_{p+q-d}(\PP_\Gamma(N))^\Gamma
\ar[d]^\cong \\
H_p(\PP_G(M))^G \otimes H_q(\PP_G(M))^G \ar[rr]^(.58)\circ & &
H_{p+q-d}(\PP_G(M))^G. }$$
\end{proof}

\section{Discrete torsion}\label{Discrete torsion}

For a gerbe over the finite group $G$, namely a 2-cocycle $\alpha \colon G \times G \to \U{1}$, one can
associate a 1-cocycle over the inertia groupoid of $G$ (see \cite{LupercioUribe6}). The inertia groupoid $\wedge G$
of $G$ is the translation groupoid of the adjoint action of $G$ on itself. Its space
of objects is $\wedge G_0 =G$, its space of arrows is $\wedge  G_1= G \times G$ and
its structural maps are $\source(g,h)=g$, $\target(g,h)=h^{-1}gh$, $\mult((g,h),(h^{-1}gh,k))= (g,hk)$.
It is Morita equivalent to the
groupoid  $\bigsqcup_{(g)} C(g)$ where $(g)$ runs over conjugacy classes of elements in $G$ and
$C(g)$ is the centralizer of $g$.

The 1-cocycle over $\wedge G$ induced by $\alpha$ is the groupoid map:
\begin{eqnarray*}
\tau \colon \wedge G & \to & \U{1} \\
(g,h) & \mapsto & \frac{\alpha(g,h)}{\alpha(h,h^{-1}gh)}
\end{eqnarray*}
where $\U{1}$ is the groupoid with one object and morphisms the unitary complex numbers.

This 1-cocycle $\tau$ once restricted to the centralizers $C(g)$ is what is known by ``discrete torsion''. For every $g \in G$ it defines  1-dimensional representation $\tau_g$ where the action of $h \in C(g)$  is given by the multiplication of $\tau(g,h)$.

The loop orbifold $\Loop \Xx:= [\PP_G(M)/G]$ is Morita equivalent to the groupoid
$$\bigsqcup_{(g)} [\PP_g(M)/C(g)].$$
As the group $G$ is finite we have that
$$H_*(\PP_G(M) \times_{G}EG;\complex) \cong H_*(\PP_G(M);\complex)^{G}$$
where the second expression means the $G$ invariant part.

\begin{definition}
The loop orbifold homology twisted by a discrete torsion $[\alpha] \in H^2(G, \U{1})$
is $$H^\alpha_*(\LL B \Xx) := \bigoplus_{(g) }\left[H_*(\PP_g(M);\complex) \otimes \tau_g \right]^{C(g)}$$
where $C(g)$ acts by translation on $H_*(\PP_g(M);\complex)$
and $\tau_g$ is the $C(g)$-representation defined above. As a vector space this is isomorphic to
$$\bigoplus_{(g) }\left[H_*(\PP_g(M);\complex) \right]^{C(g)}$$
where  for $h \in C(g)$ and $x \in H_*(\PP_g(M);\complex)$,
$h\cdot x= \tau(g,h) h_*x$.

\end{definition}

\section{Operadic structure}\label{Operadic structure}

In \cite[section 2]{CohenJones} was shown that the BV algebra
structure of the homology of free loop space of a manifold could
be understood via some suitable action of the cactus operad. In
this section we will argue that the homology of the orbifold loops
$\PP_G(M)$ can also be endowed with the action of an operad that
we have named the {\it marked $G$ chord diagram} operad. The
elements of this operad can be understood as $G$-principal bundles
over {\it marked chord diagrams}, this last operad turns out to be
homeomorphic to the cactus operad.


\subsection{Chord diagrams}
The  chord diagrams we will be interested in will have only one
boundary circle, will be flat (namely, the chords do not intersect
each other) and will have several marked points. These marked
points will be very useful when considering orbifold loops. First
we define what we mean by a {\it labeled chord diagram}.

\begin{definition}
A {\it labeled chord diagram} $c$ with $n-1$ chords consists of the following information:
\begin{itemize}
\item A counterclockwise oriented circle $S^1$ centered at the
origin in $\real^2$  with perimeter $1$. Namely $$S^1:= \{z \in
\complex :|z|=1/2\pi\}.$$ \item We will say that  $u =1 \in S^1$
is the initial marked point. \item Marked points $x_i, y_i \in
S^1$, $i \in \{1, 2, \dots, n-1\}$ with $x_i \neq y_i$, such that
they satisfy the following properties:
\begin{itemize}
\item The segment of line (known as the $i$-th {\it chord}) from
$x_i$ to $y_i$ do not intersect the $j$-th chord in the interior
of the disc (we will say in this case that the chord diagram in
flat or unknotted). We will call by $f_i \colon \{x_i\} \to
\{y_i\}$ the function that associates $y_i$ to $x_i$. We do this
because we want to think of a chord as isomorphism between two
points in the circle.  \item The $n$ connected regions that the
chords define inside the circle are numbered from $1$ to $n$.
\item The boundary of each of the $n$ connected regions that the
chords define should intersect the circle $S^1$ in a set of
measure bigger than zero.
\end{itemize}
\end{itemize}
The space of labeled chord diagrams with $n-1$ chords will be denoted by $\lcd(n)$.
\end{definition}

The term {\it labeled} stands for the labeling of the chords. The
topology of $\lcd(n)$ is the one inherited from $(S^1)^{2n-2}$.

\begin{rem} To
have a more intuitive picture of the above construction, one
could think of each chord as the place where the circle gets
``pinched" to form a cactus (following the terminology of Voronov
\cite[Section 2.2]{CohenVoronov}); each region becomes a {\it
lobe} in the cactus.
 The condition that the intersection of the
boundary of each region with the circumference of the circle must
be of measure bigger than zero ensures that each region defines a
lobe of the cactus; if the boundary of a region had only isolated
points of the circle, it would mean that the lobe got contracted
until it disappeared. We need the number of lobes to remain fixed
for each $n$ (see picture below). We would like to emphasize that
{\emph{we will not use the cactus operad in this paper,}} but
rather an equivalent operad of chord diagrams defined in this
section.
\begin{eqnarray*}
\includegraphics[height=2.5in]{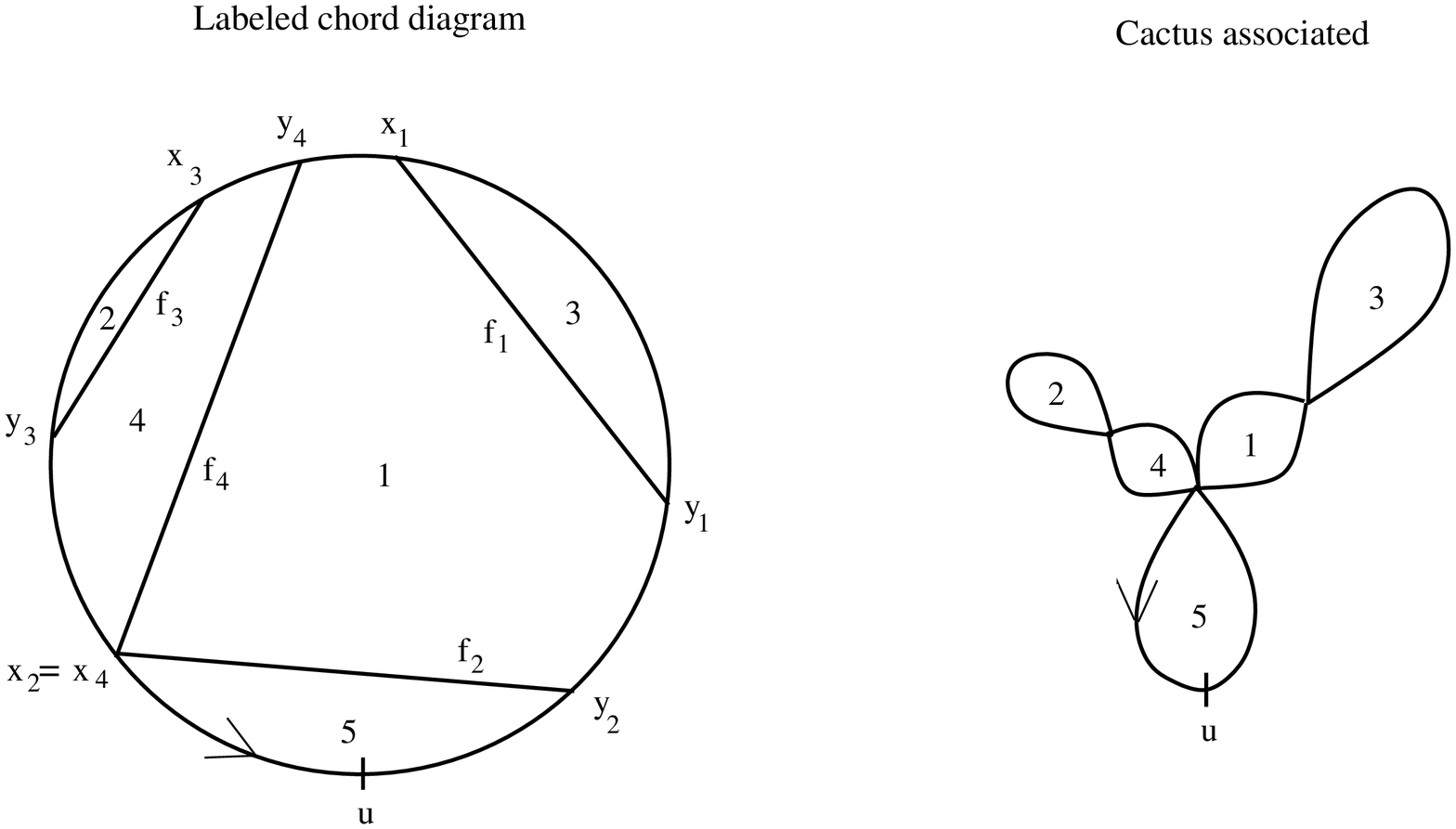} \label{labeled_cd}
\end{eqnarray*}
\end{rem}

If we erase the circle and keep the chords and the $x_i$'s,
$y_i$'s, we get a graph $\mbox{\it{graph}}(c)$, the chords are the edges and
the $x_i$ and $y_i$'s are the vertices. The condition that the
boundary of each region must contain a set of measure bigger than
zero of the circle is the same as saying that $\mbox{\it{graph}}(c)$ is a {\it
forest}. A forest is a graph whose connected components are {\it
trees}. A tree is a connected graph without loops. Note that the
points $x_i$ and $y_j$ could be equal, but they should define a
forest. All these concepts easier to understand with a picture:
\begin{eqnarray*}
\includegraphics[height=2.5in]{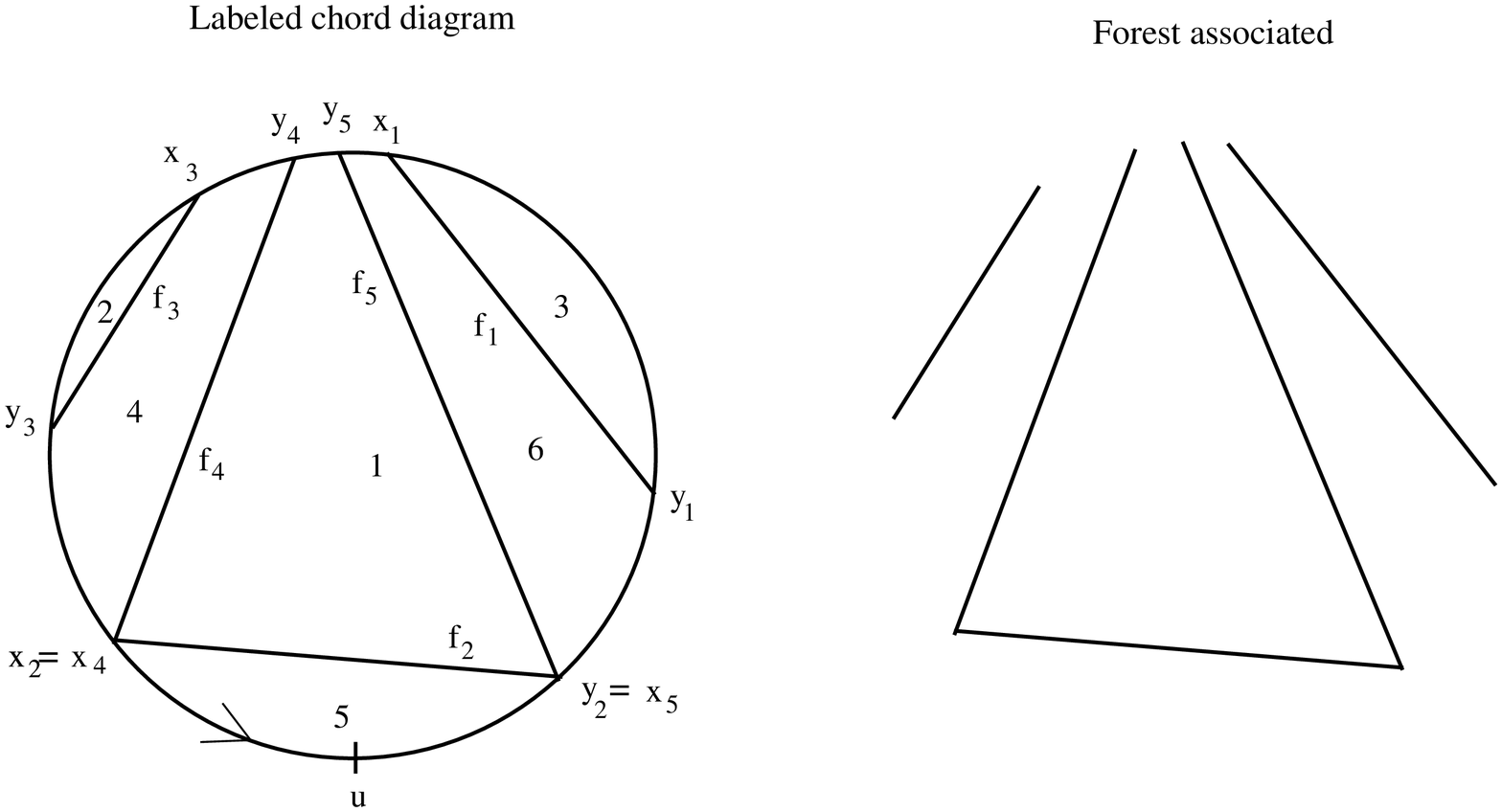} \label{forest_cd}
\end{eqnarray*}

The group $\Gamma_n \colon = \gr{S}_{n-1} \ltimes
(\gr{S}_2)^{n-1}$ acts on $\lcd(n)$ by permuting the chords, i.e.
each $\gr{S}_2$ permutes the $x_i$ with the $y_i$ and the
$\gr{S}_{n-1}$ permutes the indices. The action is free.

\begin{definition}
The space of {\it unlabeled chord diagrams } will be denoted by
$\cd(n)$ and is defined as the quotient $\cd(n) := \lcd(n) /
\Gamma_n$.
\end{definition}

Many unlabeled chord diagrams could define the same cactus as in
the following figure:

\begin{eqnarray*}
\includegraphics[height=1.8in]{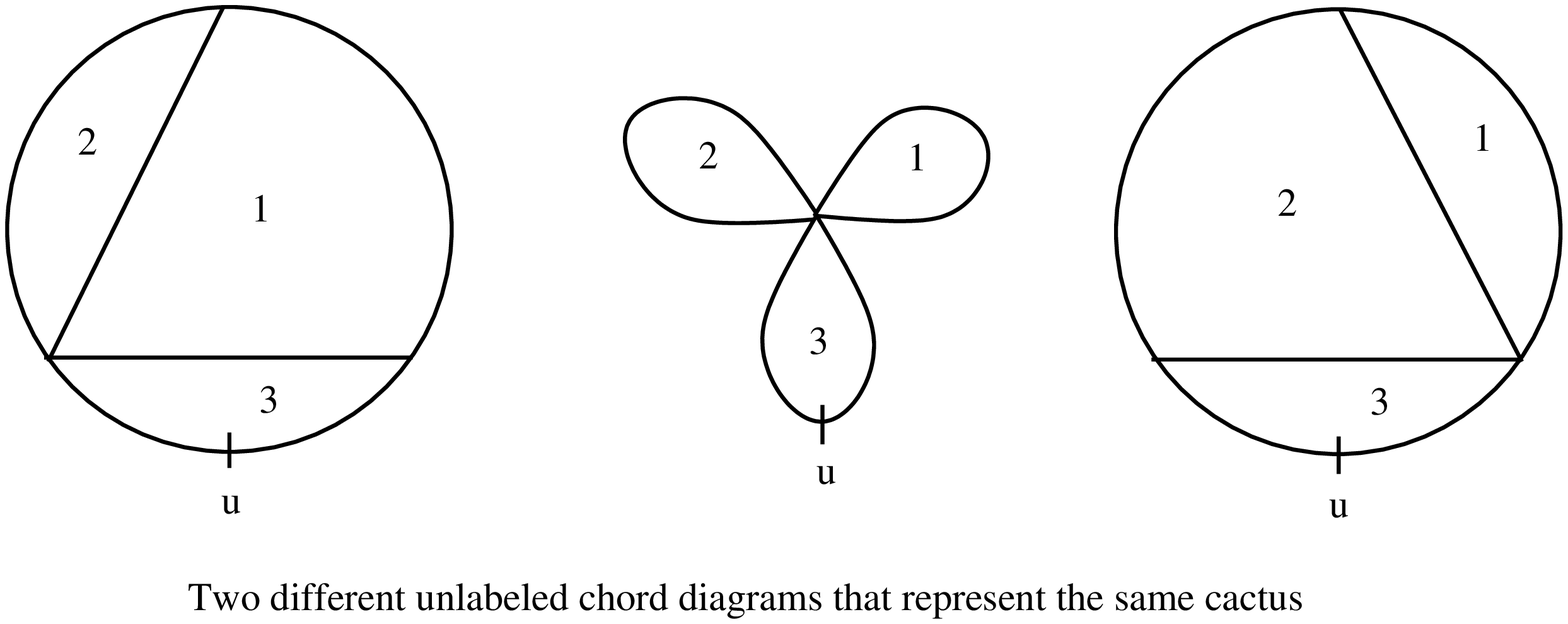} \label{equivalent_forest}
\end{eqnarray*}

For an unlabeled chord diagram $c \in \cd(n)$, the vertices of the
forest $\mbox{\it{graph}}(c)$ are divided into \emph{clusters}, namely two
vertices of $\mbox{\it{graph}}(c)$ are in the same cluster if they are in the
same tree.

The number of each region labels the different intervals in which
the circumference of the circle is divided by the chords.
Therefore one can label each such interval with the number of the
region on which it lays.

\begin{definition} Two unlabeled chord diagrams $c, c' \in \cd(n)$ have related forests if:
\begin{itemize} \item The vertices of the forests of $\mbox{\it{graph}}(c)$ and of $\mbox{\it{graph}}(c')$ are the same, and they both define the same cluster.
\item They both define the same labels for the  intervals of the boundary of the circle.
\end{itemize}

We will denote this relation by $c \sim_f c'$.

\end{definition}

We have that $\sim_f$ is  an equivalence relation. If the forest
of $\mbox{\it{graph}}(c)$ consists of $n-1$ trees with one branch, then its
equivalence class has only one element.

\begin{definition}
The space of equivalence classes of unlabeled chord diagrams $\dd
(n) := \cd(n) / \sim_f$ will be called the space of {\it chord
diagrams}.
\end{definition}

We claim that the homotopy type of $\dd(n)$ is the same as the set
of configurations of $n$ ordered  points in $\real^2$, and
therefore homotopically equivalent to the space of configurations
of $n$ little discs. Unfortunately the spaces $\{\dd(n)\}_n$ are
not endowed with an operad structure as are the little discs, for
the loops defining each region do not have a marked point.

\subsection{Marked chord diagrams}

In each of the loops determined by each of the regions of a labeled chord diagram we are going to place a marking point. This point will only be allowed to move in the loop defined by the respective region.

\begin{definition}
A {\it  marked labeled chord diagram} ${\bf c}=(c,(z_1, \dots,
z_n))$is a labeled chord diagram $c$ together with marked points
$z_i \in S^1$, $i \in \{1,\dots,n\}$, such that $z_i$ lies on the
closed interval of the circle with label $i$.
\end{definition}

As we want the $z_i$ to move only on the $i$-th loop determined by
the $i$-th region, we need to identify marked labeled chord
diagrams whenever one has that $z_i$ equals $x_j$ or $y_j$ for
some $j$. So we say that ${\bf c}$ and ${\bf c'}$ are related,
${\bf c} \sim_0 {\bf c'}$, if $c=c'$ and for all $i$, $z_i \neq
z_i'$ implies that both $z_i$ and $z_i'$ are vertices of the same
tree of the forest of $\mbox{\it{graph}}(c)$. It is clear that $\sim_0$ is an
equivalence relationship. Abusing the notation we define:

\begin{definition}
The space of {\it marked labeled chord diagrams} will be denoted by $\mlcd(n)$ and is defined as the quotient
$$\mlcd(n) := \{ {\bf c} | {\bf c}  \ \mbox{is a marked labeled chord diagram} \} / \sim_0$$
\end{definition}

\begin{eqnarray*}
\includegraphics[height=3in]{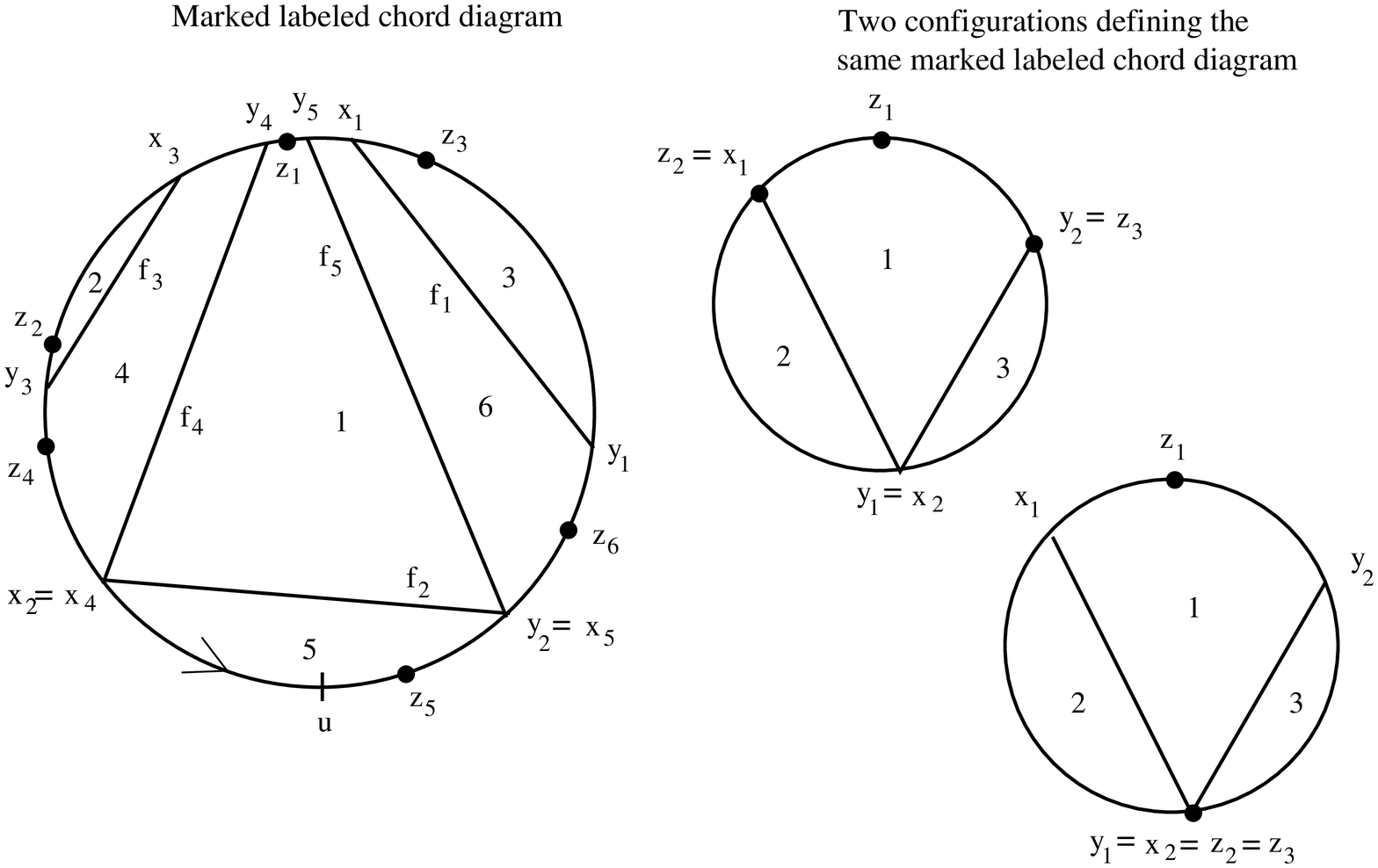} \label{marked_lcd}
\end{eqnarray*}

Then the natural forgetful map $\mu_n \colon \mlcd(n) \to
\lcd(n)$, ${\bf c} \mapsto c$ is the quotient map by the natural
action of $(S^1)^n$. Moreover, the group $\Gamma_n$ acts in the
natural way on $\mlcd(n)$ making the map $\mu_n$ into a
$\Gamma_n$-equivariant map.

\begin{definition}
The space of {\it marked unlabeled chord diagrams} will be denoted
by $\mcd(n)$ and is defined as the quotient $\mcd(n):=
\mlcd(n)/\Gamma_n$.

If furthermore we mod out  by the relation $\sim_f$ (that is well
defined in $\mcd(n)$ because of the relationship $\sim_0$) we get,
\begin{definition}
The space of {\it marked  chord diagrams} will be denoted by
$\mdd(n)$ and is defined as the quotient $\mdd(n):  = \mcd(n) /
\sim_f$.
\end{definition}
\end{definition}
Therefore we have the commutative diagram
$$\xymatrix{ \mlcd(n) \ar[d]^{\mu_n} \ar[r]_{/\Gamma_n} & \mcd(n) \ar[d]^{\mu_n'} \ar[r]_{/\sim_f} & \mdd(n) \ar[d]^{\bar{\mu}_n}\\
\lcd(n) \ar[r]_{/\Gamma_n} & \cd(n) \ar[r]_{/\sim_f} & \dd(n);
}$$
 where the forgetful maps $\mu_n' \colon \mcd(n) \to \cd(n)$ and $\bar{\mu}_n \colon \mdd(n) \to \dd(n)$ are quotient maps by the action of $(S^1)^n$. This is true because the action of the group $\Gamma_n$ is free in both spaces, and because the relation $\sim_f$ in $\mcd(n)$ is defined by the lift of the relation $\sim_f$ in $\cd(n)$.

\vspace{.5cm}

We will argue now that each marked unlabeled chord diagram
determines uniquely a cactus (with marked points). So, let us
recall the definition of the cactus operad $Cac$ given in
\cite[Section 2]{CohenJones}. A point $a$ in the space $Cac(n)$ is
a collection of $n$ oriented, parameterized circles $a_1, \dots,
a_n$ with radii $r_i$ so that $\sum_{i=1}^n r_i =1$. Each circle
has a marked point $s_i \in a_i$ given by the image under the
parametrization of the base point $1 \in S^1$. The circles can
intersect each other at a finite number of points creating a
``cactus type configuration", namely that the dual graph of this
configuration is a tree (the dual graph of the configuration has
for vertices the circles and an edge between two circles whenever
they intersect, we will denote it by $tree(a)$). The boundary of
the cactus (the union of the circles) is equipped with a basepoint
$w$ together with a choice of which component the base point $w$
lies in. The edges coming into any vertex are also equipped with a
cyclic ordering. The topology of the space of cacti $Cac(n)$ is
described in \cite{CohenVoronov, Voronov}.

\begin{proposition}  \label{proposition map delta}
There is a surjective map  $\delta_n \colon \mcd(n) \to Cac(n)$ that factors through
$$\xymatrix{ \mcd(n) \ar[rr]^{\delta_n} \ar[rd]_{\bar{\mu}_n} && Cac(n) \\
& \mdd(n) \ar[ru]_{\bar{\delta}_n} &
}$$
Moreover, the map $\bar{\delta}_n \colon \mdd(n) \to Cac(n)$ is a homeomorphism.
\end{proposition}

\begin{proof}
A  marked unlabeled chord diagram ${\bf c}=(c,(z_1, \dots, z_n))$
has $n$ connected regions in its interior, each one with a number.
Once the chords of the boundary of each region are contracted, one
obtains a loop $\alpha_i$  for each region, with certain perimeter
$r_i$. Moreover, each loop $\alpha_i$ has a marked point $z_i$.
For each $i$, take $a_i$ to be the parameterized circle of radius
$r_i$ and marked point $s_i \in a_i$ (the image of $u=1\in S^1$ in
the parametrization), and construct the unique orientation
preserving map $h_i \colon \alpha_i \to a_i$ of constant velocity
and winding number $1$ , such that $h_i(z_i)=s_i$. As the total
perimeter of the loops $\alpha_i$ is $1$, then
$\sum_{i=1}^nr_i=1$.

The intersection of the circles $a_i$ are obtained via  the maps
$h_i$ and the points where the chords got contracted. The cyclic
order of the circles at each intersection point is given by the
following algorithm. Intersection points in the cactus correspond
to trees in $\mbox{\it{graph}}(c)$. Fix a tree. Each vertex of a given tree
has a small neighborhood in $S^1$ that has two labels. These
labels in turn correspond to labels in the lobes of the cactus
that intersect at the point corresponding to the tree. We will
order the labels via the counterclockwise orientation of $S^1$ at
the vertex of the tree. In this manner the lobes get a cyclic
ordering.

Finally, let $j$ be the region where $u=1\in S^1$ lies on. In the
case that $u$ is over a vertex of the forest $\mbox{\it{graph}}(c)$, take $j$
to be the region where $u$ would go if it followed the
orientation. Then take $w := h_j(u)$ lying in the circle $a_j$.

The cactus $a$ is uniquely determined by the marked labeled chord
diagram ${\bf c}$, then we define $\delta_n({\bf c}) := a$.

If two marked unlabeled chord diagrams ${\bf c}$ and ${\bf c'}$
satisfy ${\bf c} \sim_f {\bf c'}$, then they define the same
cactus, i.e. $\delta_n({\bf c}) = \delta_n({\bf c'})$. For the
intersection points of the cactus are defined by each of the trees
of the forests. When the trees collapse to a point, both ${\bf c}$
and ${\bf c'}$ define the same intersection points. Then we can
define $\bar{\delta}_n([{\bf c}]) := \delta_n({\bf c})$.

To show that $\bar{\delta}_n$ is a homeomorphism we will construct
its inverse map. Take a cactus $a$ and let $S^1$ be the
parameterized circle of perimeter $1$ with marked point $u$.
Construct the unique orientation preserving map $\kappa \colon a
\to S^1$ with constant velocity $\frac{1}{2\pi}$  that maps the
boundary of the cactus $a$ to the circle $S^1$ such that
$\kappa(w)=u$. This map is unique because we start in $w$
following the boundary of the cactus in the direction of the
orientation. If $w$ lies in an intersection point, since we know
in which lobe $w$ lies, we move through the boundary of that lobe.
Once an intersection point is reached, we use the cyclic order to
know where to continue. The image under $\kappa$ of the points
$s_i$ define the numbering of the regions as well as the marked
points $z_i$. The image under $\kappa$ of each one of the
intersection points define the vertices of each of the trees of
the forest. This information defines uniquely an element in
$\mdd(n)$. Therefore we have defined a continuous map $\sigma_n
\colon Cac(n) \to \mdd(n)$. It is straightforward to check that
$\sigma_n$ is the inverse of $\bar{\delta}_n$.
\end{proof}

Notice that if there are $n-1$ intersection points in the cactus
$a$, then we have $n-1$ pairs of points in the circles $S^1$ and
there is only one way we can fill them with trees. In this case we
can find a marked unlabeled chord diagram for the cactus $a$. If
there are less than $n-1$ intersection points, there is no
canonical way to decide how to construct the trees of the forest
of the marked unlabeled chord diagram. Nevertheless there is
always a tree, take for example the linear graph defined by the
cyclic order given by the cyclic orientation starting at any point
in a cluster.

\vspace{0.5cm}

 The spaces $\{\mlcd(n)\}_{n \geq 0} $, $\{\mcd(n)\}_{n \geq 0} $ and $\{\mdd(n)\}_{n \geq 0} $ form {\it operads}. An operad (see \cite[Def. 2.3]{Voronov}) is a set of spaces $\{\OO(n)\}_n$ with  the following properties:
 \begin{itemize}
 \item A composition law:
 $$\OO(m) \times \OO(k_1) \times \cdots \times \OO(k_m) \To \OO(k_1 + \cdots + k_m)$$
 \item A right action of the symmetric group $\gr{S}_n$ on $\OO(n)$.
 \item A unit $e \in \OO(1)$.
 \end{itemize}
 such that the composition is associative and equivariant with respect to the symmetric groups, and
 the identity satisfies natural properties with respect to the composition.

 Let's see how  the operad structure is defined in $\mlcd$. The action of $\gr{S}_n$ is given by permutations of the subindices of the marked points $z_i$ (and therefore the numbering of the regions), and the composition is obtained by
adding the chords of the diagrams ${\bf{c}_i} \in \mlcd(k_i)$ into
the diagram of ${\bf{c}} \in \mlcd(m)$. This is done as follows:
the diagram $\bf{c}_i$ will be {\it patched} into the $i$-th
region (where the point $z_i$ lies in $\bf{c}$). $\bf{c}_i$ is
linearly rescaled to match the length of the $i^{th}$ loop of
$\bf{c}$. The patching starts by sending the marked point $u^i$ of
$\bf{c}_i$ to the point $z_i$ and one continues using the
orientations. The points $\{x^i_j,y^i_j\}_{1 \leq j < k_i}$ are
patched then onto $\bf{c}$ thus defining an element $\td{\bf{c}}
\in \mlcd(k_1 + \cdots + k_m)$. The numbering of the chords and of
the regions in $\td{\bf{c}}$ is obtained by starting with the
chords and the regions in $\bf{c}_1$, then the ones in $\bf{c}_2$
and so forth. The marking $u$ of $\bf{c}$ remains the same in
$\td{\bf{c}}$. The identity $e \in \mlcd(1)$ is the chord diagram
with $z_1 =u$.

The operad structure of $\mlcd$ induces  the operad structures on
$\mcd$ and $\mdd$.

\begin{theorem}
The cactus operad $Cac$ and the marked  chord diagram operad
$\mdd$ are isomorphic via the maps $\bar{\delta}=
\{\bar{\delta}_n\}_n$.
\end{theorem}

\begin{proof}
In proposition \ref{proposition map delta} we showed that the maps $\bar{\delta}_n \colon \mlcd(n) \to Cac(n)$ are homeomorphisms. The fact that they $\gr{S}_n$ equivariant and compatible with the operad structure is straightforward.
\end{proof}

In \cite{CohenGodin, CohenJones} it  is shown how the homology of
the Cacti operad  can act on the homology of the free loop space
of a manifold thus realizing the BV-algebra structure on the
homology of the loops defined by Chas and Sullivan
\cite{ChasSullivan}. This uses the result of Voronov
\cite{Voronov} that shows that the cactus operad is equivalent to
the framed little disc operad, and the result of Getzler
\cite{Getzler} that shows that the algebraic structure of a
BV-algebra is captured by an action of the framed little disc
operad.

\subsection{Marked G-chord diagrams}

 Now we want to construct the moduli space of principal $G$ bundles over the marked chord diagrams,
 i.e. principal $G$ bundles over $\bf{c}$ for all $\bf{c}$'s. We need
 to lift all the markings of the marked labeled chord diagram in order to have a  manifold and not an orbifold.
 This fact can be easily seen on the moduli space of principal $G$ bundles over a circle. If one marks
 a lift of $0$ then the isomorphism class  can be characterized by the holonomy, en element in $G$, but without the marking
 the only invariant we can measure is the conjugacy class of the holonomy. The first one is a manifold, the second one is
 an orbifold.

First, let's fix for each $g \in G$ a $G$-principal bundle $\pi
\colon Q_g \to S^1$ (as in section \ref{principal bundles} with a
marked point $\tilde{u} \in Q_g$, over the parameterized circle of
total arclength $1$, such that $\pi$ is a local isometry,
$\pi(\tilde{u})=u$ and the holonomy of $Q_g$ starting at
$\tilde{u}$ is $g$.

 \begin{definition}
 A {\it $G$ marked labeled chord diagram} ${\bf W}$ with holonomy $g$ over the marked labeled chord diagram ${\bf c} \in \mlcd(n)$ consist of the following information:
 \begin{itemize}
 \item The principal bundle $\pi \colon Q_g \to S^1$.
 \item The marked labeled chord diagram ${\bf c}$.
 \item Lifts $ \td{z}_1, \dots, \td{z}_n $ in $Q_g$ for the marked points $ z_1, \dots, z_n$ in ${\bf c}$.
 \item Lifts $\td{x}_i, \td{y}_i \in Q_g$, $1 \leq i <n$ of $x_i$ and $y_i$ in ${\bf c}$ , together with isomorphisms
 $\td{f}_i \colon \pi^{-1}(x_i) \stackrel{\cong}{\to} \pi^{-1}(y_i)$,
 $1 \leq   i <n$ that are $G$-equivariant and such that $\td{f}_i(\td{x}_i) = \td{y}_i$
 \end{itemize}
 We will denote the space of such $G$ marked labeled chord diagrams $${\bf W}=(Q_g, {\bf c}, (\td{x}_i)_i, (\td{y}_i)_i,(\td{z}_i)_i )$$ by $\gmlcd(n,g)$.

Considering all the possible holonomies we define $\gmlcd(n) :=
\bigsqcup_{g \in G} \gmlcd(n,g).$
 \end{definition}

\begin{example}
 In the graph below we can see a marked $\gr{S}_3$ chord diagram over
  the marked labeled chord diagram of one chord. The dashed double arrows
   represent the isomorphism $\td{f}_1$ and the pointed vertical
 lines represent the fibers of the points $x_1$ and $y_1$. We have identified
 the fiber of $x_1$ with the elements
 of $\gr{S}_3$ that appear on the left of the figure, and the fiber of $y_1$
  with the elements on the right.

\begin{eqnarray}
\includegraphics[height=4.0in]{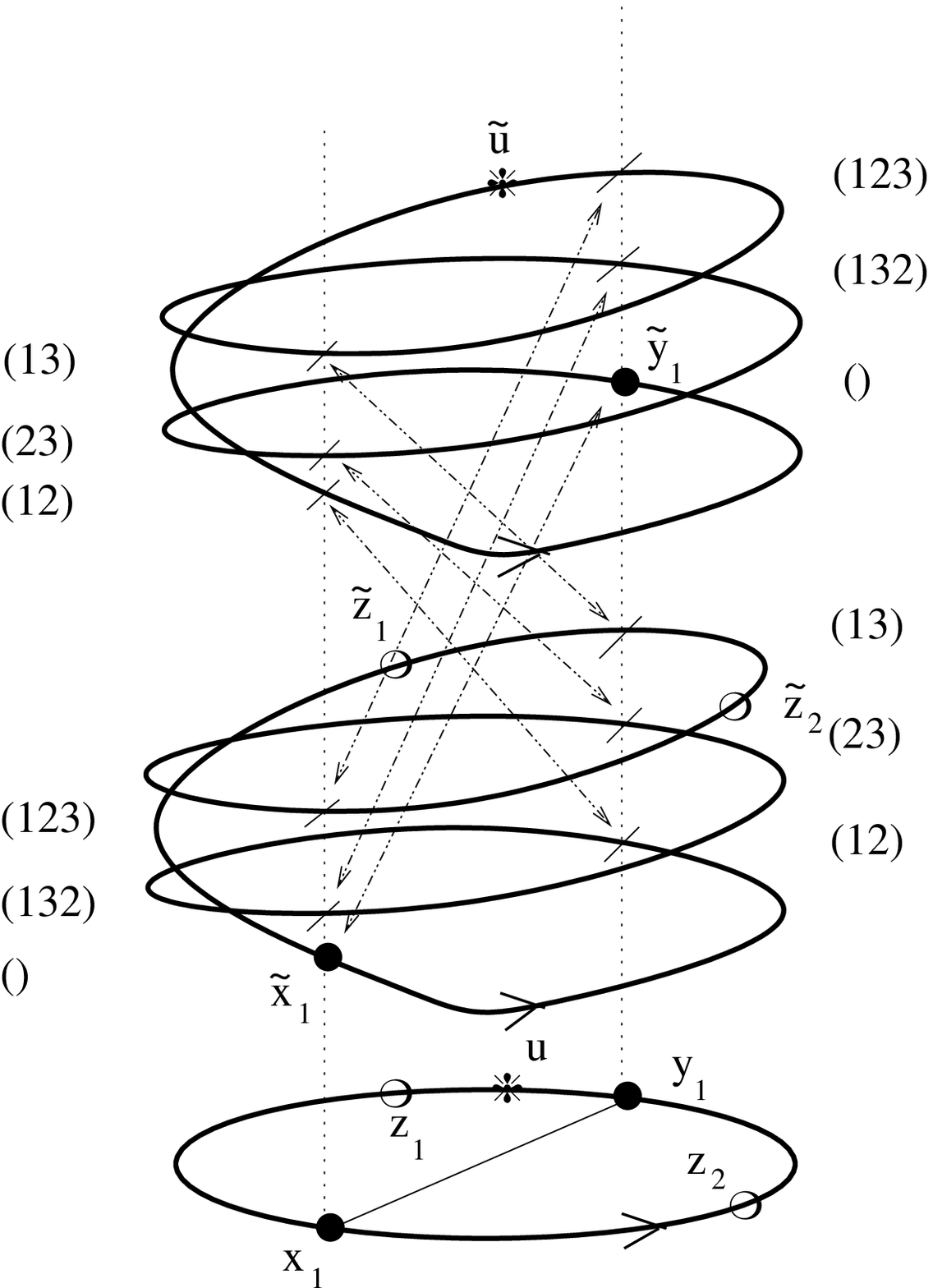} \label{gbundle}
\end{eqnarray}

From the figure one can see that the holonomy of the outer circle
(forgetting the chords) starting at $\td{u}$  is $(1,3,2)$, the
holonomy of the region numbered 1 starting at $\td{z}_1$ is
$(2,3)$ and of the region 2 starting at $\td{z}_2$ is also
$(2,3)$. It is worth pointing out here that if one takes the
$\gr{S}_3$ principal bundle associated to either one of the
regions, each one of them would have 3 connected components. Each
component will give a different holonomy, namely $(1,2)$, $(2,3)$
or $(1,3)$; all the conjugates of $(2,3)$. All these holonomies
can be attained by changing the position of the points $\td{z}_i$.
\end{example}

\begin{proposition}
The projection maps $\gmlcd(n,g) \to \mlcd(n)$, ${\bf
W}\mapsto {\bf c}$, are principal $ G^{3n-2}$ bundles.
\end{proposition}
\begin{proof}
There are $|G|$ choices for each of the $\tilde{x}_i,
\tilde{y}_i$, and $ \tilde{z}_i$'s. The element  $$( (g_i)_i,
(h_i)_i,(k_i)_i) \in G^{n-1} \times G^{n-1} \times G^n $$ sends
${\bf W}=(Q_g, {\bf c}, (\td{x}_i)_i, (\td{y}_i)_i,(\td{z}_i)_i )$
to ${\bf W'}=(Q_g, {\bf c}, (\td{x}_i g_i)_i, (\td{y}_i
h_i)_i,(\td{z}_i k_i )_i )$.
\end{proof}

Therefore for all $g \in G$, $$\gmlcd(n,g)/G^{3n-2} = \mlcd(n).$$

\vspace{0.5cm}

We want to consider now principal bundles over the marked
unlabeled chord diagrams and we will use $\gmcd(n)$ to denote this
moduli space.

Consider the diagonal action
\begin{equation} \label{diagonal action on gmlcd}
\gmlcd(n,g) \times G^{n-1}  \to \gmlcd(n,g)
\end{equation}
$$\ \ \ \ \ \  (Q_g, {\bf c}, (\td{x}_i)_i, (\td{y}_i)_i,(\td{z}_i)_i ) \times
(g_i)_i  \mapsto  (Q_g, {\bf c}, (\td{x}_i g_i)_i, (\td{y}_i
g_i)_i,(\td{z}_i )_i )
$$

\begin{definition}
Let $\bargmlcd(n,g) := \gmlcd(n,g)/G^{n-1}$ where the  $G^{n-1}$
action is the diagonal action of \ref{diagonal action on gmlcd}.
\end{definition}
 An equivalence class ${\bf \overline{W}}$ in $\bargmlcd(n,g)$ consist of  of a $G$ marked labeled chord diagram where we keep track of the $G$ equivariant maps $\tilde{f}_i \colon \pi^{-1}(\{x_i\}) \to \pi^{-1}(\{y_i\})$ lifting the $f_i$'s and the $\tilde{z}_i$, but there are no marked points on the fibers of the $x_i, y_i$'s.

If $G^{n-1}$ acts on $\bargmlcd(n,g)$ by sending the maps
$(\tilde{f}_i)_i$ to the maps $(\tilde{f}'_i)_i$, where
$\tilde{f}'_i(v) := \tilde{f}_i(v)g_i$, one gets that
$\bargmlcd(n,g)$ is a $G^{2n-1}$ principal bundle over $\mlcd(n)$.

Over the space $\bargmlcd(n,g)$ the group $\Gamma_n$ acts by
permuting the $\tilde{f}_i$'s lifting the action of $\Gamma_n$ on
$\mlcd(n)$.

\begin{definition}
The space of {\it $G$ marked unlabeled chord diagrams} with
holonomy $g$,  will be denoted by $\gmcd(n,g)$ and is defined as
the quotient $$\gmcd(n,g) := \bargmlcd(n,g) / \Gamma_n.$$
\end{definition}
The following lemma follows immediately from the definition.
\begin{lemma} \label{lemma forgetful map in gmcd}
The forgetful map $\gmcd(n,g) \to \mcd(n)$ that sends the $G$ marked unlabeled chord diagram ${\bf W}$ to its base (a marked unlabeled chord diagram ${\bf c}$ ) is finite map of order $|G|^{2n-1}$.
\end{lemma}

Notice that there is no canonical way to get an action of the
group $G^{n-1}$ on $\gmcd(n,g)$; nevertheless $G^n$ acts on
$\gmcd(n,g)$ permuting the $\tilde{z}_i$'s.

\vspace{0.5cm}

An element ${\bf W}$ in $\gmcd(n,g)$ consists of a principal
$G$-bundle over the marked unlabeled chord diagram ${\bf c}$,
together with $G$-equivariant maps  $\pi^{-1}(a) \to \pi^{-1}(b)$,
where  $a,b$ are vertices of the same chord in ${\bf c}$. The
forest of $\mbox{\it{graph}}({\bf c})$ divides the vertices of the chords of
${\bf c}$ into \emph{clusters}. The inverse image under $\pi$ in
${\bf W}$ of each of these clusters is further divided into $|G|$
\emph{subclusters} using the functions $\{\tilde{f}_i\}_i$. We
define then

\begin{definition}
Two $G$ marked unlabeled chord diagrams ${\bf W}, {\bf W'} \in \gmcd(n,g)$ are forest related ${\bf W} \sim_F {\bf W'}$, if:
\begin{itemize}
\item their marked unlabeled chord diagrams ${\bf c}, {\bf c'}$
are forest related, i.e. ${\bf c} \sim_f {\bf c'}$, \item the
subclusters of the inverse image under $\pi$ of the vertices of
$\mbox{\it{graph}}({\bf c})$ are the same as the connected components of the
inverse image under $\pi'$, and \item $\tilde{z}_i=\tilde{z}'_i$
for all $i$ (recall that the $\tilde{z}_i$'s as well as the
$\tilde{z}'_i$'s belong to $Q_g$).
\end{itemize}
\end{definition}

By definition $\sim_F$ is an equivalence relation, so we define
\begin{definition} \label{definition gmdd}
The space of {\it $G$ marked chord diagrams} with holonomy $g$,
will be denoted by $\gmdd(n,g)$ and is defined as the quotient
$$\gmdd(n,g) := \gmcd(n,g) / \sim_F.$$
\end{definition}
From lemma \ref{lemma forgetful map in gmcd} and definition \ref{definition gmdd} we get
\begin{lemma}
The forgetful map $\gmdd(n,g) \to \mdd(n)$ that sends the $G$ marked chord diagram ${\bf W}$ to its base (a marked chord diagram ${\bf c}$ ) is finite map of order $|G|^{2n-1}$.
\end{lemma}

Therefore we get the commutativity of the following diagram
$$\xymatrix{
\gmlcd(n,g) \ar[r]_{/G^{n-1}} \ar[rd] & \bargmlcd(n,g) \ar[r]_{/\Gamma_n} \ar[d] & \gmcd(n,g) \ar[r]_{/ \sim_F} \ar[d] & \gmdd(n,g) \ar[d] \\
& \mlcd(n) \ar[r]_{/\Gamma_n} & \mcd(n) \ar[r]_{\sim_f} & \mdd(n)
}$$ where all the vertical maps are finite or order $|G|^{2n-1}$,
theleft diagonal map is a $G^{3n-2}$ fibration, and the top
horizontal maps are $G^n$ equivariant (by acting of the
$\tilde{z}_i$'s).

We can now measure the holonomy around each of the regions that defines a marked chord diagram. Then let
 \begin{eqnarray*} \ih \colon \gmlcd(n,g)  & \to & G^n\\
\BF{W} &\mapsto & (h_1, h_2, \dots , h_n)
\end{eqnarray*}
be the map that assigns the holonomy around each of the  $n$ regions. This is measured  by
 starting from the points $\td{z}_i$ following the
induced orientation. Whenever a point $\td{x}_j$ is reached the
path is continued on $\td{y}_j$ (or viceversa) and so forth until
one reaches a point on the same fiber as $\td{z}_i$. The name
$\ih$ stands for \emph{incoming holonom}. We will denote by
\begin{eqnarray*} \oh \colon \gmlcd(n)  & \to & G\\
\BF{W} &\mapsto & g
\end{eqnarray*}
the map that assigns $g$ to ${\bf W}$ whenever ${\bf W} $ belongs
to $\gmlcd(n,g)$. Recall that this is  the total holonomy of
$\BF{W}$ measured starting from $\td{u}$; the name $\oh$ stands
for \emph{outgoing holonomy}.

The incoming and outgoing holonomies are defined similarly on each
of the spaces $\bargmlcd, \gmcd, \gmdd$. Then we can define
\begin{definition}
For $g, h_1, \dots h_n \in G$ with $\mathbf{h} = ( h_1, \dots,
h_n)$ let
$$\gmlcd(n,g,\mathbf{h}) := \{ \BF{W} \in \gmlcd(n) \ |
 \ \ih(\BF{W}) = (h_1, \dots, h_n) \ {\rm and } \ \oh(\BF{W}) =g \}.$$
 Similarly for $\bargmlcd, \gmcd, \gmdd$.
\end{definition}

So, $\gmlcd(n,g) = \bigsqcup_{\mathbf{h} } \gmlcd(n,g,\mathbf{h})$
and similarly for $\bargmlcd, \gmcd, \gmdd$. Notice that some of
the sets $\gmlcd(n,g,\mathbf{h})$ may be empty. This is because
the inner holonomy $\mathbf{h}$ plus the geometric information of
the bundle ${\bf W}$ determines uniquely the outgoing holonomy.

Fixing the incoming and outgoing holonomies we get the following
commutative diagram \begin{eqnarray} \label{commutative diagrams with (n,g,h)}
\xymatrix{
  & \bargmlcd(n,g, \mathbf{h}) \ar[r]_{/\Gamma_n} \ar[d]^{/G^n} & \gmcd(n,g, \mathbf{h}) \ar[r]_{/ \sim_F} \ar[d]^{/G^n} & \gmdd(n,g, \mathbf{h}) \ar[d]^{/G^n} \\
\gmlcd(n,g, \mathbf{h}) \ar[ru]_{/G^{n-1}} \ar[r]  & \mlcd(n) \ar[r]_{/\Gamma_n} & \mcd(n) \ar[r]_{\sim_f} & \mdd(n)
}
\end{eqnarray} where all the vertical maps are $G^n$ principal
bundles (the action is on the $\tilde{z}_i$'s).

Summarizing for the spaces $\gmdd$ we have:
\begin{lemma} \label{lemma equivalences  gmdd}
$$\gmdd(n) = \bigsqcup_{(g, \mathbf{h})} \gmdd(n,g,\mathbf{h})$$
and $$\gmdd(n,g,\mathbf{h})/G^n = \mdd(n).$$ \end{lemma}

\vspace{0.5cm}

The \emph{$G$-graded operad action} for $\gmdd=\{\gmdd(n)\}_n$
given by maps:

\begin{eqnarray}
\gmdd(n) \  {}_{\ih} \! \! \times_{\{\prod \oh\}} \prod_{j=1}^n
\gmdd(l_j) \To \gmdd(l_1 + l_2 + \cdots + l_n) \label{G
operadic structure}
\end{eqnarray}
 is defined in
the natural way such that it is compatible with the one of
$\{\mdd(n)\}_n$, where the set $\gmdd(n) \  {}_{\ih} \! \!
\times_{\{\prod \oh\}} \prod_{j=1}^n \gmdd(l_j)$ is the one that
makes the following a cartesian square
$$
\xymatrix{ \gmdd(n) \  {}_{\ih} \! \! \times_{\{\prod \oh\}}
\prod_{j=1}^n \gmdd(l_j) \ar[rr] \ar[d] & &  \prod_{j=1}^n
\gmdd(l_j) \ar[d]^{\prod \oh} \\
 \gmdd(n) \ar[rr]^{\ih} & & G^n
 },$$
and where the identity $\tilde{e} \in \gmdd(1)$ is the trivial $G$
principal bundle over $e\in \mdd(1)$ with $\tilde{z}_i =
\tilde{u}$.

\vspace{0.5cm}

We have decided to call the structure described in Lemma\
\ref{lemma equivalences  gmdd} \ a \emph{$G$-graded operad}. All
the operads in this paper with a $G$ in their name are $G$-graded
operads.

Now we are ready to construct the action of this operad in the
homology of the loop orbifold.

Recall that the loop orbifold can be seen as $[\Bun_G(S^1,M)/G]$
where $$\Bun_G(S^1,M) = \bigsqcup_{g \in G} \Bun_g(S^1,M)$$ and
$\Bun_g(S^1,M)$ is the set of $G$-equivariant maps from $Q_g$ to
$M$ with $Q_g$ the $G$-principal bundle over the circle with
holonomy $g$ (measured from $q_0$ the lift of $0$) and the
projection being a local isometry.

\begin{definition}
For $\BF{W} \in \gmlcd(n,g,\mathbf{h})$, let $L_{\bf{W}}M$ be
$$L_{\bf{W}}M := \{ \gamma \colon Q_g\to M \ | \ \gamma \
{\rm is } \ G-{\rm equivariant} \ {\rm and} \ \gamma(\td{x}_i) =
\gamma(\td{y}_i) \ 1 \leq i < n \}.$$
\end{definition}
Notice that the space $L_{\bf{W}}M$ could also have been defined
as the space of $G$ equivariant maps $\gamma \colon Q_g \to M$
such that $\gamma(\tilde{a})=\gamma(\tilde{b})$ where  $\tilde{a},
\tilde{b}$ are any two points in $Q_g$ lifting vertices $a,b \in
graph(c)$ and such that $\tilde{a}, \tilde{b}$ are in the same
subcluster of $\pi^{-1}(graph(c))$ defined via the $f_i$'s.
Therefore we have:

\begin{lemma} \label{lemma same space L}
If ${\BF{W}} , {\bf{W'}}\in \gmlcd(n,g,\mathbf{h})$ are $G$ marked
labeled chord diagram that define the same $G$ marked chord
diagram $\bf{Z}$ (i.e. ${\BF{W}} \mapsto {\bf{Z}}$ and  ${\bf{W'}}
\mapsto \bf{Z}$ under the map $\gmlcd(n,g, \mathbf{h}) \to
\gmdd(g,n,\mathbf{h})$ of (\ref{commutative diagrams with
(n,g,h)})), then $L_{\bf{W}}M$ is equal to $L_{\bf{W'}}M$.
\end{lemma}

 By forgetting the information about the marked points (except $\tilde{u}$),
  we get the natural map
$$\theta_{\BF{W}} \colon L_{\bf{W}}M \to \Bun_g(S^1,M).$$

The space $L_{\bf{W}}M$ can also be viewed as the pullback of an
evaluating mapping of the product $\prod_{j} \Bun_{h_j}(S^1,M)$
defined as follows. The marked $G$-chord diagram $\BF{W}$ induces
$G$-principal bundles $W_i$ over the circles $a_i$ of radius
$\frac{1}{2\pi}$ of the perimeter of the $i$-th region, by taking
the $G$ principal bundle over the $i$-th region that $\bf W$
defines by contracting its chords via the maps $f_j$. The marked
point being $\td{z}_i$.

 On each of these bundles $W_i$ denote by $\alpha_l$ the points on $W_i$ that
correspond to $\td{x}_l$ (or $\td{y}_l$) where $l\in I_i$ with
$I_i$ the set of $j$ such that $\td{x}_j$ (or $\td{y}_j$) is on
the $i$-th region. Let $m_i$ be the number of elements in $I_i$,
and as all the chords touch two regions we have that $m_1+ \cdots
+ m_n = 2(n-1)$.

Define the evaluation map
$$ev_{\BF{W}} \colon \prod_{i=1}^n \Bun_{h_i}(S^1,M) \To
(M)^{2(n-1)}$$ as follows. Let $s_i\colon S^1 \to a_i$ be the
identification of the unit circle with $a_i$ obtained by scaling
down the unit circle so as to have the radius of $a_i$, and
rotating it so the base point $0 \in S^1$ is mapped to the marked
point $z_i \in a_i$. And take $\td{s}_i \colon Q_{h_i} \to W_i$ to
be the identification of the corresponding $G$-principal bundles,
compatible with $s_i$ and where $q_0$ is mapped to $\td{z}_i$. Let
$\beta_l, l \in I_i$ the points on $Q_{h_i}$ corresponding to
$\alpha_l, l \in I_i$ under the map $\td{s}_i$. Define
\begin{eqnarray*}
ev_{W_i} \colon \Bun_{h_i}(S^1,M) & \to & (M)^{m_i}\\
\sigma \colon Q_{h_i} \to M & \mapsto & (\sigma(\beta_{j^i_1}),
\dots , \sigma(\beta_{j^i_{m_i}}))
\end{eqnarray*}
with $j^i_1 < j^i_2 < \cdots < j^i_{m_i}$ all of them in $I_i$.

Now define
$$ev_{\BF{W}} = ev_{W_1} \times \cdots \times ev_{W_n} \colon
\prod_{i=1}^{n} \Bun_{h_i}(S^1,M) \to (M)^{m_1}\times \cdots
\times (M)^{m_n} = (M)^{2(n-1)}.$$

As each of the chords touch two regions, then each of the $j$'s $1
\leq j <n$ belong to two of the $I_i$'s. This induces a diagonal
map
$$\Delta_{\BF{W}} \colon (M)^{(n-1)} \to (M)^{m_1}\times \cdots
\times (M)^{m_n} = (M)^{2(n-1)}$$ defined by
$$\Delta_{\BF{W}}(p_1, \dots , p_{n-1}) = (q_{j^1_1},\dots,
q_{j^1_{m_1}}, q_{j^2_1} , \dots, q_{j^2_{m_2}}, \dots
q_{j^{n}_1},\dots, q_{j^{n}_{m_{n}}})$$ with $p_j= q_{j^i_k}$
whenever $j =j^i_k$. As noted above, for each $j$ there are only
two $j^i_k$'s that are equal to $j$.

Note that the map $\Delta_{\BF{W}}$ is the same as the $(n-1)$-th
product of the diagonal map $M \to M\times M$ composed with a
permutation of the group $\gr{S}_{2(n-1)}$.

Now observe that the following is a cartesian pullback square:

$$\xymatrix{
L_{\bf{W}}M \ar[d]_{ev_{int}} \ar[rr]^{\td{\Delta}_{\BF{W}}} & &
\prod_{i=1}^n \Bun_{h_i}(S^1,M) \ar[d]^{ev_{\BF{W}}} \\
(M)^{(n-1)} \ar[rr]_{\Delta_{\BF{W}}} && (M)^{2(n-1)} }
$$
where $ev_{int} \colon L_{\BF{W}}M \to (M)^{n-1} $ evaluates a map
$\gamma \colon {\BF{W}} \to M$ at the $n-1$ vertices of the
chords, namely $ev_{int}(\gamma) = (\gamma(\td{x}_1), \dots,
\gamma(\td{x}_{n-1}))$.

The normal bundle $\eta(\Delta_{\BF{W}})$ of the diagonal
embedding $\Delta_{\BF{W}}$ is then isomorphic to $TM^{n-1} \to
(M)^{n-1}$. By the tubular neighborhood theorem, we have proven
the following:

\begin{lemma} \label{lemma tubular neighborhood}
The image of the embedding $$\td{\Delta}_{\BF{W}} \colon
L_{\BF{W}}M \to \prod_{i=1}^n \Bun_{h_i}(S^1,M)$$ has an open
neighborhood homeomorphic to the total space of the pullback
$ev_{int}^*TM^{n-1}$
\end{lemma}

We now consider the above construction for families of diagrams
by letting $\BF{W}$ vary in $\gmlcd(n,g,\mathbf{h})$. Consider the
set
$$L_{(n,g,\mathbf{h})}M = \{ (\BF{W}, \gamma) \colon \BF{W} \in
\gmlcd(n,g,\mathbf{h}) \ {\rm and} \ \gamma \in L_{\BF{W}}M \},$$
the map
\begin{eqnarray*}
\td{\Delta}_{(n,g, \mathbf{h})} \colon L_{(n,g,\mathbf{h})}M & \to
&
\gmlcd(n,g, \mathbf{h}) \times \prod_{i=1}^n \Bun_{h_i}(S^1,M) \\
(\BF{W}, \gamma) & \mapsto & (\BF{W},
\td{\Delta}_{\BF{W}}(\gamma)),
\end{eqnarray*}
the map
\begin{eqnarray*}
ev_{(n,g, \mathbf{h})} \colon L_{(n,g,\mathbf{h})}M & \to &
M^{n-1} \\
(\BF{W}, \gamma) & \mapsto & ev_{int}(\gamma) =
(\gamma(\tilde{x}_1), \dots ,\gamma(\tilde{x}_{n-1})),
\end{eqnarray*}
and the map
\begin{eqnarray} \label{definition theta}
\theta_{(n,g,\mathbf{h})} \colon L_{(n,g,\mathbf{h})}M & \to &
\Bun_g(S^1,M) \\
(\BF{W}, \gamma) & \mapsto & \theta_{\BF{W}}(\gamma).
\end{eqnarray}
Defining $\xi(n,g,\mathbf{h}) := ev_{(n,g,\mathbf{h})}^*TM^{n-1}$,
the $d(n-1)$ dimensional vector bundle over $L_{(n,g,\mathbf{h})}$
whose fiber over $(\BF{W},\gamma)$ is the sum of tangent spaces
$$\xi(n,g,\mathbf{h})|_{(\BF{W},\gamma)} =
\bigoplus_{i=1}^{n-1}T_{\gamma(\td{x}_i)}M,$$ together with lemma
\ref{lemma tubular neighborhood}, we have that

\begin{lemma} \label{lemma embedding gmlcd}
The image of the embedding $$\td{\Delta}_{(n,g, \mathbf{h})}
\colon L_{(n,g,\mathbf{h})}M \to \gmlcd(n,g, \mathbf{h}) \times
\prod_{i=1}^n \Bun_{h_i}(S^1,M)$$ has an open neighborhood
homeomorphic to the total space $\xi(n,g,\mathbf{h})$.
\end{lemma}

\vspace{0.5cm}

Now we need to define a similar map as in lemma \ref{lemma
embedding gmlcd} but for the space $\gmdd$. This we will do by
stages defining it first for $\bargmlcd$ and then on $\gmcd$.

Recall that $$\bargmlcd(n,g,\mathbf{h}) =
\gmlcd(n,g,\mathbf{h})/G^{n-1}.$$ Take ${\bf W}, {\bf W'} \in
\gmlcd(n,g,\mathbf{h})$ such that both map to ${\bf V} \in
\bargmlcd(n,g,\mathbf{h})$. By lemma \ref{lemma same space L} we
have that $L_{\bf W}M = L_{\bf W'}M$, therefore we can define an
action of $G^{n-1}$ on $L_{(n,g,\mathbf{h})}M$ that takes $({\bf
W}, \gamma) \mapsto ({\bf W'}, \gamma)$ making the embedding
$\td{\Delta}_{(n,g, \mathbf{h})}$ $G^{n-1}$ equivariant. This
action can be lifted to the vector bundle $\xi(n,g,\mathbf{h})$ by
using the isomorphism
$$\xi(n,g,\mathbf{h})|_{(\BF{W},\gamma)} =
\bigoplus_{i=1}^{n-1}T_{\gamma(\td{x}_i)}M \cong
\xi(n,g,\mathbf{h})|_{(\BF{W'},\gamma)} =
\bigoplus_{i=1}^{n-1}T_{\gamma(\td{x}'_i)}M.$$ Therefore we have
that $\xi(n,g,\mathbf{h})$ becomes a $G^{n-1}$ equivariant vector
bundle over $L_{(n,g,\mathbf{h})}M$.

\begin{definition}
Let $\bar{\xi}(n,g,\mathbf{h}) := {\xi}(n,g,\mathbf{h})/G^{n-1}$
be the vector bundle over $\bar{L}_{(n,g,\mathbf{h})}M :=
{L}_{(n,g,\mathbf{h})}M/G^{n-1}$.
\end{definition}
We have then

\begin{lemma}
The induced map
$$\bar{\Delta}_{(n,g, \mathbf{h})}
\colon \bar{L}_{(n,g,\mathbf{h})}M \to \bargmlcd(n,g, \mathbf{h})
\times \prod_{i=1}^n \Bun_{h_i}(S^1,M)$$ is an embedding and has
an open neighborhood homeomorphic to the total space
$\bar{\xi}(n,g,\mathbf{h})$.
\end{lemma}

Notice that for the vector bundle $\bar{\xi}(n,g,\mathbf{h})$, its
fibers are
$$\bar{\xi}(n,g,\mathbf{h})|_{(\BF{V},\gamma)} \cong \left(
\bigoplus_{(g_1,\dots,g_{n-1}) \in G^{n-1}}
\bigoplus_{i=1}^{n-1}T_{\gamma(\td{x}_i g_i)}M \right)^{G^{n-1}}
$$

\vspace{0.5cm}

Recall that the group $\Gamma_n$ acts freely on
$\bargmlcd(n,g,\mathbf{h})$ and that its quotient defines
$\gmcd(n,g,\mathbf{h})$. For $\sigma \in \Gamma_n$, the action of
$\sigma$ on $\bar{L}_{(n,g,\mathbf{h})}M$ takes $({\bf V},\gamma)$
to $({\bf V}\sigma, \gamma)$ and this action can be lifted to an
action on $\bar{\xi}(n,g,\mathbf{h})$ where the canonical
isomorphism
$$\bar{\xi}(n,g,\mathbf{h})|_{(\BF{V},\gamma)} \cong
\bar{\xi}(n,g,\mathbf{h})|_{(\BF{V}\sigma,\gamma)}$$ is
$$\left(
\bigoplus_{(g_1,\dots,g_{n-1}) \in G^{n-1}}
\bigoplus_{i=1}^{n-1}T_{\gamma(\td{x}_i g_i)}M \right)^{G^{n-1}}
\cong \left( \bigoplus_{(g_1,\dots,g_{n-1}) \in G^{n-1}}
\bigoplus_{i=1}^{n-1}T_{\gamma((\td{x}_i g_i)\sigma)}M
\right)^{G^{n-1}}.$$

We have then that the vector bundle $\bar{\xi}(n,g,\mathbf{h})$ is
$\Gamma_n$ equivariant.

\begin{definition}
Let $\hat{\xi}(n,g,\mathbf{h}) :=
\bar{\xi}(n,g,\mathbf{h})/\Gamma_n$ be the vector bundle over
$$\hat{L}_{(n,g,\mathbf{h})}M :=
\bar{L}_{(n,g,\mathbf{h})}M/\Gamma_n.$$
\end{definition}
We have then

\begin{lemma}
The induced map
$$\hat{\Delta}_{(n,g, \mathbf{h})}
\colon \hat{L}_{(n,g,\mathbf{h})}M \to \gmcd(n,g, \mathbf{h})
\times \prod_{i=1}^n \Bun_{h_i}(S^1,M)$$ is an embedding and has
an open neighborhood homeomorphic to the total space
$\hat{\xi}(n,g,\mathbf{h})$.
\end{lemma}

Notice that for the vector bundle $\hat{\xi}(n,g,\mathbf{h})$, its
fibers are
$$\hat{\xi}(n,g,\mathbf{h})|_{(\BF{Y},\gamma)} \cong \left[
\bigoplus_{\sigma \in \Gamma_n} \left(
\bigoplus_{(g_1,\dots,g_{n-1}) \in G^{n-1}}
\bigoplus_{i=1}^{n-1}T_{\gamma((\td{x}_i g_i)\sigma)}M
\right)^{G^{n-1}} \right]^{\Gamma_n}
$$
where ${\bf V}$ maps to ${\bf Y}$ under the map
$\bargmlcd(n,g,\mathbf{h}) \to \gmcd(n,g,\mathbf{h})$.

\vspace{0.5cm}  {
 The equivalence relation $\sim_F$ is easily
defined on $\hat{L}_{(n,g,\mathbf{h})}M$. For $({\bf Y},\gamma),
({\bf Y'},\gamma') \in \hat{L}_{(n,g,\mathbf{h})}M$, we say that
$({\bf Y},\gamma) \sim_F ({\bf Y'},\gamma')$ if and only if ${\bf
Y} \sim_F {\bf Y'}$ and $\gamma = \gamma'$. We now want to lift
this relation to $\hat{\xi}(n,g,\mathbf{h})$. Since the fibers of
$\hat{\xi}(n,g,\mathbf{h})$ over $({\bf Y},\gamma)$ and $({\bf
Y'},\gamma)$ are canonically isomorphic  whenever ${\bf Y} \sim_F
{\bf Y'}$ (this is because the images of $\gamma$ under ${\bf Y}$
and ${\bf Y'}$ are the same) then this isomorphism defines the
equivalence relation $\sim_F$ on $\hat{\xi}(n,g,\mathbf{h})$. }

\begin{definition}
Let $\check{\xi}(n,g,\mathbf{h}) :=
\hat{\xi}(n,g,\mathbf{h})/\sim_F$ be the vector bundle over
$$\check{L}_{(n,g,\mathbf{h})}M :=
\hat{L}_{(n,g,\mathbf{h})}M/\sim_F.$$ And let
$$\check{\theta}_{(n,g,\mathbf{h})} \colon
\check{L}_{(n,g,\mathbf{h})}M\to \Bun_g(S^1,M)$$ be the
composition of the maps
$$\xymatrix{{L}_{(n,g,\mathbf{h})}M \ar[r] & \check{L}_{(n,g,\mathbf{h})}M
\ar[rr]^{\theta_{(n,g,\mathbf{h})}} &&  \Bun_g(S^1,M)}$$where
$\theta_{(n,g,\mathbf{h})}$ is defined in (\ref{definition
theta}).
\end{definition}
We can conclude,

\begin{lemma}
The induced map
$$\check{\Delta}_{(n,g, \mathbf{h})}
\colon \check{L}_{(n,g,\mathbf{h})}M \to \gmdd(n,g, \mathbf{h})
\times \prod_{i=1}^n \Bun_{h_i}(S^1,M)$$ is an embedding and has
an open neighborhood homeomorphic to the total space
$\check{\xi}(n,g,\mathbf{h})$.
\end{lemma}

Notice that this allows us to perform a Pontrjagin-Thom collapse
map
$$ \tau \colon \gmdd(n,g, \mathbf{h}) \times
\prod_{i=1}^n \Bun_{h_i}(S^1,M) \to (\check{L}_{(n,g,\mathbf{h})}M
)^{\check{\xi}(n,\mathbf{h})}$$ that induces a homomorphism in
homology
$$H_p\left(\gmdd(n,g, \mathbf{h}) \times
\prod_{i=1}^n \Bun_{h_i}(S^1,M)\right)\to H_{p-d(n-1)}(
\check{L}_{(n,g,\mathbf{h})}M)$$  that once composed with the
homomorphism
$$
\xymatrix{ H_*(\check{L}_{(n,g,\mathbf{h})}M)
\ar[rr]^{(\check{\theta}_{(n,g,\mathbf{h})})_*}&&
H_*(\Bun_h(S^1,M))}$$ produces the following.

\begin{proposition} \label{G operad structure map}
There is a homomorphism in homology
$$H_p\left(\gmdd(n,g, \mathbf{h}) \times
\prod_{i=1}^n \Bun_{h_i}(S^1,M)\right) \to
H_{p-d(n-1)}(\Bun_g(S^1,M))$$ which induces an action of the
homology of the $G$-graded operad $\gmdd$.
\end{proposition}

\begin{proof}
The compatibility is clear from the operadic structure of $\gmdd$
(see \ref{G operadic structure}).
\end{proof}

Now we are ready to show that the homology with rational
coefficients of $\LL B \Xx= \LL M_G$ has a BV algebra structure,
this we will prove by using the proposition \ref{G operad
structure map} together with the equivalences of lemma \ref{lemma
equivalences  gmdd}.

\begin{theorem}
The homology with rational coefficients of the free loop space of
the Borel construction $[M/G]$, i.e. $H_*(\LL M_G , \rational)$,
has a Batalyn-Vilkovisky algebra structure. Moreover when the
group is the identity group coincides with the Chas-Sullivan
construction.
\end{theorem}

\begin{proof}
Let's start by recalling that $\LL M_G$ is homotopically
equivalent to the space $\Bun_G(S^1,M) \times_G EG$, and  as $G$ is finite
we have that
$$H_*(\LL M_G , \rational) \cong H_*( \Bun_G(S^1,M) \times_G EG, \rational)
\cong H_*( \Bun_G(S^1,M) , \rational)^G$$ where the last group
means the $G$-invariant part of $H_*( \Bun_G(S^1,M) , \rational)$.

Using corollary \ref{lemma equivalences  gmdd} we have the
following set of maps:

$$H_*\left(\mchord(n) \times \prod_{i=1}^n (\LL M_G), \rational \right)
\cong  H_*\left(\mchord(n) \times \prod_{i=1}^n (\Bun_G(S^1,M)/G)
, \rational\right) \To $$
$$ \frac{1}{|G|^n} \bigoplus_{g\in G, \mathbf{h} \in G^n}
H_*\left(\gmdd(n,g, \mathbf{h})/G^{n} \times \prod_{i=1}^n
\left(\Bun_G(S^1,M)/G\right) , \rational\right) \cong $$
$$ \frac{1}{|G|^n}
\bigoplus_{g\in G, \mathbf{h} \in G^n} H_*\left(\gmdd(n,g,
\mathbf{h}), \rational\right)^{G^{n}} \otimes \bigotimes_{i=1}^n
H_* \left(\Bun_G(S^1,M) ,\rational\right)^G \To $$
$$  \frac{1}{|G|^n} \bigoplus_{g\in G, \mathbf{h} \in G^n}
 H_*\left(\gmdd(n,g, \mathbf{h}), \rational\right) \otimes
\bigotimes_{i=1}^n H_* \left(\Bun_G(S^1,M) ,\rational\right)^G \To
$$
$$
H_* \left(\Bun_G(S^1,M) ,\rational\right)^G$$ where the last
homomorphism is the one from proposition \ref{G operad structure
map} whenever is defined; namely, when the inner holonomy of ${\bf
W}$ matches the holonomy of the loops in $\prod_{i=1}^n
\Bun_G(S^1,M)$ then the product is the one from proposition \ref{G
operad structure map}, otherwise is zero.

The only thing left to prove is that the image is indeed
$G$-invariant, but this should be clear from the fact that we have
added over all possible $G$ marked chord diagrams on the
definition of the homomorphism.

Therefore we have the set of  homomorphisms

$$H_p\left(\mdd(n) \times \prod_{i=1}^n (\Bun_G(S^1,M)/G)
, \rational\right) \To H_{p-d(n-1)}\left(\Bun_G(S^1,M)/G ,
\rational \right),$$

that are compatible with the operadic structure of $\mdd$.

By proposition \ref{proposition map delta} the chord diagram
operad is isomorphic to the cactus operad, which by a theorem of
Voronov \cite[Theorem 2.3]{Voronov} is homotopy equivalent to the
framed little disk operad, we have that $H_*(\LL M_G, \rational)$
becomes a Batalyn-Vilkovisky algebra.

When the group $G$ is the trivial group, then
$\Bun_G(S^1,M) = \LL M$ and $\gmdd = \mdd$, and we recover the
homological action of the cactus operad in the homology of $\LL M$
defined in \cite{CohenJones}.

\end{proof}

\begin{rem}
Note that the action of $\mdd(1)$ on $$\Bun_G(S^1,M)/G =
\PP_G(M)/G = \coarse{\Loop \Xx},$$ is the one induced by the $S^1$
action defined on \ref{circle action}.
\end{rem}

\section{Computations}\label{Computations}

\begin{example}
Let $M$ be a smooth manifold and consider $\Xx=[M/\{1\}]$ (in other words we consider the case when $G=\{1\}$). Then it is clear that $\PP_g(M) = \PP_G(M) = \LL M$ is simple the free loop space and $H_*(\Loop \Xx) = H_*(\LL M).$
By the work of Cohen and Jones we recover the Chas-Sullivan BV-algebra in this case.
\end{example}

\begin{example}
Let $G$ be a finite group and consider $\Xx=[\bullet/G]$ be the
orbifold consisting of a point $M= \bullet$ being acted by $G$.
Sometimes this orbifold is denoted by $\BB G$ (not to be confused with
$B G$ the classifying space of $G$).
Clearly every loop and every path in this case is constant, namely the
space $\PP_g(M) = \star_g$ is a point, and so $\PP_G(M)$ is in
one-to-one correspondence with $G$. Therefore the category
$[\PP_G(M)/G]$ is equivalent to the category $[G/G]$ of $G$ acting on
$G$ by conjugation, for we have $$h(\star_g)=\star_{hgh^{-1}}.$$ For
each $g\in G$ the stabilizer of this action is the centralizer
$$C(g)=\{ h\in G | hgh^{-1}=g \}.$$ Now, in the category $[G/G]$ an
object $g\in G$ is isomorphic to $g'\in G$ if and only if $g$ and $g'$
are conjugate. Therefore we have the equivalence of categories $$\Loop \Xx \simeq [\PP_G(M)/G]\simeq [G/G]
\simeq \coprod_{(g)} [\star_g/C(g)].$$ Here $(g)$ runs through the
conjugacy classes of elements in $g \in G$.
From this we can conclude that the equivalence $$\LL B \Xx = B \Loop \Xx $$ becomes in this particular case (cf. \cite{LupercioUribe5})
$$\LL B G \simeq \coprod_{(g)} B C(g)$$
This equation becomes at the level of homology with complex coefficients the center of the group algebra $$H_*(\LL B G) \cong Z(\complex[G])$$
and in fact $H_*(\LL B G)$ is simply the Frobenius algebra of Dijkgraaf and Witten \cite{DijkgraafWitten}.

The reader may be interested in comparing this result with that of \cite{AbbaspourCohenGruher}.

\end{example}

Let $X$ be a topological space endowed with the action of a
connected Lie group $\Gamma$. Take $G \subset \Gamma$ finite and
consider the quotient $X/G$ and the map $\pi : X \to X/G$ the
projection.

\begin{lemma} \label{lemma Lie group1}
The projection map induces an isomorphism
$$\pi_* : H_*(X ; \rational) \stackrel{\cong}{\to} H_*(X/G; \rational).$$
\end{lemma}

 \begin{proof}
Take $g \in G$ and its induced action $g: X \to X$. We claim that
$g_* : H_*(X) \stackrel{=}{\to} H_*(X)$ is the identity. Join the
identity of $\Gamma$ with $g$ with a path $\alpha_t \in \Gamma$
(i.e. $\alpha_0= id_\Gamma$ and $ \alpha_1= g$), hence $\alpha_t$
is a homotopy  between the identity and $g$, therefore $g_* = id$.

Taking the averaging operator

\begin{eqnarray*}
H_*(X;\rational) & \stackrel{\alpha}{\to} & H_*(X;\rational)^{G} \\
x & \mapsto& \left(\frac{1}{|G|} \sum_{g \in G} g_* x \right) (=
x)
\end{eqnarray*}
and using that
 $ H_*(X;\rational)^{G} \stackrel{\iota}{\cong} H_*(X/G; \rational) $
the isomorphism follows, for it is not hard to check that $\pi_*=\iota \alpha$.
\end{proof}

With the same hypothesis as before consider now the orbifold
loops, namely $\PP_g X = \{ f : [0,1] \to X | f(0)g = f(1) \}$.

\begin{lemma} \label{lemma Lie group2}
There is a $C(g)$-equivariant homotopy equivalence between $\LL X$
and $\PP_g X$.
\end{lemma}

\begin{proof}

Let $\alpha_t : [0,1] \to G$ be the map defined in lemma \ref{lemma Lie group1}. Consider the maps
\begin{eqnarray}
\rho : \PP_g X  \to \LL X  \ \ \ \ \ \ \mbox{and} \ \ \ \ \
\tau : \LL X \to  \PP_g X \label{homotopy tau}
\end{eqnarray}
where
\begin{eqnarray*}
\rho(f)(s):= \left\{
\begin{array}{ccc}
f(2s) & \mbox{if} & 0 \leq s \leq \frac{1}{2} \\
f(1) \alpha_{2s-1}^{-1} & \mbox{if} & \frac{1}{2} \leq s \leq 1
\end{array}
\right.
\end{eqnarray*}
and
\begin{eqnarray*}
\tau(\sigma)(s):= \left\{
\begin{array}{ccc}
\sigma(2s) & \mbox{if} & 0 \leq s \leq \frac{1}{2} \\
\sigma(1) \alpha_{2s-1} & \mbox{if} & \frac{1}{2} \leq s \leq 1.
\end{array}
\right.
\end{eqnarray*}

The composition  $ \rho  \circ \tau : \LL X \to \LL X$ is clearly  homotopic to the identity. The same holds for $\tau
\circ \rho$. The maps $\rho$ and $\tau$ are trivially $C(g)$-equivariant.

\end{proof}

\begin{cor} \label{corollary Lie group}
The group structure of the loop homology of $[X/G]$ can be seen as
$$H_*(\Loop [X/G]; \rational) \cong
 \bigoplus_{(g)} H_*(\LL X ;\rational).$$
\end{cor}
\begin{proof}
It follows from the lemmas \ref{lemma Lie group1} and \ref{lemma Lie group2} and the fact that
$$H_*(\Loop [X/G]; \rational) \cong \bigoplus_{(g)} H_*(\PP_g X / C(g);\rational).$$

\end{proof}
{\bf Notation:} Let $\Xx$ be an orbifold of dimension $d$. Let's denote the loop homology of $\Xx$ by
$$\quaternion_*(\Xx) := H_{*+d}( \Loop \Xx).$$
In this way the orbifold string product $\quaternion_*(\Xx)$ is graded associative.

\begin{example} \label{example Lens}
The loop homology of the lens spaces $L_{(n,p)} = S^n / \integer_p$ ($n$ odd, $p>0$) is
$$\quaternion_*(L_{(n,p)}) =H_*(\LL L_{(n,p)}) = \Lambda[a] \otimes \rational[u,v]/(v^p=1)$$
with $a \in \quaternion_{-n}(L_{(n,p)})$, $v \in \quaternion_{0}(L_{(n,p)})$ and $u \in \quaternion_{n-1}(L_{(n,p)})$.
\end{example}
\begin{proof}
As the action of $\integer_p$ on $S^n$ comes from the action of $S^1$ on $S^n$ via the Hopf fibration, we can use
corollary \label{corollary Lie group} . Let $g $ be a generator of $\integer_p$, then
$$\quaternion_*(L_{(n,p)}) \cong \quaternion_*([S^n / \integer_p]) \cong \bigoplus_{j=0}^{p-1} H_*(\PP_{g^j} S^n)^{\integer_p},$$
as graded vector spaces.

As $H_*(\PP_{g^j} S^n)^{\integer_p} \cong H_*(\PP_{g^j} S^n)$ the string product $\circ$ could be calculated from the following commutative diagram
$$\xymatrix{
H_*(\PP_{g^j} S^n)^{\integer_p} \times H_*(\PP_{g^k} S^n)^{\integer_p} \ar[d]_\cong
\ar[rr]^(.6)\circ & & H_*(\PP_{g^{j+k} }S^n)^{\integer_p} \ar[d]^\cong \\
H_*(\PP_{g^j} S^n) \times H_*(\PP_{g^k} S^n) \ar[rr]^(.6)\circ && H_*(\PP_{g^{j+k}} S^n).
}$$

The map $\tau^j : \LL S^n \to \PP_{g^j} S^n$ defined in
(\ref{homotopy tau}) gives an isomorphism in homology, so we can
define the generators of the homology of $\PP_{g^j} S^n$ via the
map $\tau^j$ and the loop homology of the sphere, namely
$\quaternion_*(S^n) = H_*(\LL S^n) \cong \Lambda[a]\otimes \rational[u]$
(see \cite{CohenJonesYan}). Denote then by $\sigma^j_k$ the
generator of the group $H_{k+n}(\PP_{g^j} S^n)$ and using that
$(\tau^j)_*$ is an isomorphism one gets that $\sigma^j_{(n-1)l -n} = \tau^j_*(a u^l)$,
$\sigma^j_{(n-1)l } =\tau^j_*(u^l)$ and $\sigma^j_m =0$ for all
other values of $m$.

We claim now that $$\sigma^j_l \circ \sigma^k_m = \sigma^{j+k}_{l+m}.$$
The identity follows from the fact that $$\sigma^j_l \circ \sigma^k_m = \tau^j_*(\sigma^0_l )\circ \tau^k_*(\sigma^0_m )
= \tau^{j+k}_*(\sigma^0_l \circ \sigma^0_m) = \tau^{j+k}_*(\sigma^0_{l+m})= \sigma^{j+k}_{l+m}$$
where the second identity follows from the definition of the maps $\tau$ and the third identity follows from
the algebraic structure of $\Lambda[a]\otimes\rational[u]$.

From this we can deduce that the map $\tau_{j*}: H_*(\LL S^n) \to H_*(\PP_{g^j} S^n)$ maps $\sigma_k^0 \mapsto \sigma_k^0 \circ \sigma_0^j$ where $\sigma_0^j$ is the $n$-simplex of paths that to every $x$ in $S^n$ assigns the path that goes from $x$ to $xg^j$ through the $S^1$ action.

We are only left to prove that when $j+k=p$ the formula $ \sigma^j_l \circ \sigma^k_m = \sigma^{0}_{l+m}$ holds.
So, let $\beta : S^n \to \LL S^n$ be the map that to a point $x$ in the sphere associates the free loop defined that starts and ends in $x$ and travels in the direction of the
$S^1$ action. Now define the map $\phi : \LL S^n \to \LL S^n$
that takes a loop $\gamma$ to $\gamma \circ \beta$. The map $\phi$ is homotopic to the identity because
the cycle $\beta$ is homotopic to the cycle of constant loops over the sphere (one way to prove this uses the fact that the odd dimensional spheres have two orthogonal never vanishing vector fields). Therefore we have that $\tau^p : \LL S^n \to \PP_{g^p}S^n = \LL S^n $ is homotopic to the identity.

We can conclude then that the elements $a= \sigma^0_{-n}$, $v=  \sigma^1_0$ and $u= \sigma^0_{n-1}$ generate the loop homology of $L_{(n,p})$, and the only extra condition is that $v^p=1$. Therefore

$$\quaternion_*(L_{(n,p)}) = \Lambda[a] \otimes \rational[u,v]/(v^p=1)$$
\end{proof}


\begin{example}
Take the orbifold defined by the action of $\integer_p$ onto $S^2$
given by rotation of $2\pi/p$ radians with respect to the
$z$-axis. Then the loop homology of $[S^2/\integer_p]$ is
$$\quaternion_*([S^2/\integer_p]) = \Lambda[b] \otimes
\rational[a,v,y]/(a^2, ab , av , y^p-1) $$
\end{example}
\begin{proof}
The action of $\integer_p$ comes from the $S^1$ action on $S^2$
given by rotation about the $z$-axis. therefore the calculation of the loop
homology product follows the same argument as in the
example \ref{example Lens}. To make the notation simpler we will work with $p=2$ (
$\integer_2 = \{1,g\}$); the other cases are similar.

From \cite{CohenJonesYan} we know that the loop homology of $S^2$
is given by
\begin{eqnarray} \quaternion_*(S^2)=\Lambda[b] \otimes \integer[a,v]/(a^2,ab,2av) \label{loop homology sphere} \end{eqnarray}
with $|b|=1$, $|a|=-2$, $|v|=2$. Since $\tau : \LL S^2 \to \PP_g
S^2$ is a homotopy equivalence, we will follow the argument of example \ref{example Lens}. The only different argument is on the behavior
of the map $\phi:= \tau^2 : \LL S^2 \to \LL S^2$. In homology,  $\phi_*$ maps $\alpha \in H_k(\LL S^2)$ to $\alpha \circ \beta \in H_k(\LL S^2)$ where $\beta \in H_2(\LL S^2) = \quaternion_0(S^2)$ is the class of the map $S^2 \to \LL S^2$ that assigns to every point $x$ the loop that starts at $x$ and rotates around the $z$ axis, and $\circ$ is the homology string product.

We claim that $\beta = 1 + av$ in the notation of (\ref{loop  homology sphere}), (the proof of this fact will be postponed to lemma \ref{lemma beta=av}). As $av$ is a torsion class, i.e. $2av=0$, then in rational homology $\phi_*$ is the identity map. As in example \ref{example Lens}, we can add a new variable $y$ that behaves like a root of unity, and we conclude that
$$\quaternion_*([S^2/\integer_2]) = \Lambda[b] \otimes
\rational[a,v,y]/(a^2, ab , av , y^2-1).$$
\end{proof}

\begin{lemma} \label{lemma beta=av}
The homology class $\beta \in H_2 (\LL S^2)=\quaternion_0(S^2)$ of the map $S^2 \to \LL S^2$  that to a point $x$ assigns the loop that starts at $x$  and winds around the sphere once by the $S^1$ action, and the homology class $1 + av \in H_2(\LL S^2)=\quaternion_0(S^2)$ as in (\ref{loop homology sphere}), are equal. \end{lemma}

\begin{proof}
When we contract all the loops of $\beta$ through the north pole we end up with the homology class $[S^2] + \xi$, where $[S^2]$ is the fundamental class of the sphere (constant loops) and therefore the unit in $1=[S^2] \in \quaternion_0(S^2)$, and $\xi$ is defined in what follows. For $\theta \in S^1$ and $P_S$ the south pole, consider the map $f :S^1 \times S^1 \to \LL S^2$ such that the function $f_\theta=f(\cdot, \theta) : S^1 \to \LL S^2$ is the loop of based loops that starts at the constant loop in $P_S$ and goes around the sphere (as a rubber band) at the angle $\theta$. The class $f_{\theta *}([S^1])$ is the generator of $H_1(\LL S^2)$, and the class $f_*([S^1\times S^1])$ is $\xi$.
We claim that $\xi = av$.

We know  that the homology spectral sequence of the Serre fibration $\Omega S^2 \to \LL S^2 \to S^2$ has for $E_2$-term
$$E_2^{p,q} = H_p(S^2) \otimes H_q(\Omega S^2)$$
with non trivial differential $d^2(u \otimes x^{2k+1})=2\iota \otimes x^{2k+2}$ where $x \in H_1(\Omega S^2)$, $\iota \in H_0(S^2)$, $1_\Omega \in H_0(\Omega S^2)$ and $u \in H_2(S^2)$ are generators respectively. Also we know from \cite{CohenJonesYan} that $av= \iota \otimes x^2$.

Denote by $\dot{T} S^2 \stackrel{\pi}{\longrightarrow} S^2$ the sphere bundle of the tangent bundle $T S^2 \to S^2$. The map $\pi$ is an $S^1$-fibration and a point in $\dot{T} S^2$ consists of a pair $(z,v)$ where $z \in S^2$ and $v$ is a unit vector tangent to $S^2$ at $z$. For each point $(z, v)$ we can define a map $h_{(z,v)} : S^1 \to \LL S^2$ in the same way that the function $f_\theta$ was defined two paragraphs above; namely, $h_{(z,v)}$ is the loop of loops that starts with the constant loop at $z$ and sweeps the sphere as a rubber band, following the direction of the oriented maximum circle tangent to the vector $v$. We can assemble all the functions $h_{(z,v)}$ by letting $(z,v)$ vary and we can obtain a function $$\psi : S^1 \times \dot{T} S^2 \to \LL S^2$$ such that $\psi(\phi, (z,v)) = h_{(z,v)}(\phi)$.

The map $\psi$ defines a map of Serre fibrations
\begin{eqnarray} \label{map of Serre fibrations}
\xymatrix{
S^1 \times S^1 \ar[r] \ar[d] & \Omega S^2 \ar[d] \\
S^1 \times \dot{T} S^2 \ar[r]^\psi \ar[d] & \LL S^2 \ar[d] \\
S^2 \ar[r]^= & S^2}\end{eqnarray}
that induces a map in spectral sequences. If $\epsilon \in H_0(S^1) \otimes H_0(S^1)$, $a \in H_1(S^1) \otimes H_0(S^1)$, $b \in H_0(S^1) \otimes H_1(S^1)$, $ c \in H_1(S^1) \otimes H_1(S^1)$, are the generators in homology, at the second term of the map of spectral sequences $$\psi_* : H_p(S^2) \otimes H_q(S^1 \times S^1) \to H_p(S^2) \otimes H_q(\Omega S^2)$$ induces the following identities:
\begin{itemize}
\item $\psi_*(\epsilon)=1_\Omega$,
\item $\psi_*(b) = 0$ and
\item $\psi_*(a) = x$ because the functions $f_\theta$ determine the generator $x$ of $H_1(\Omega S^2)$.
\end{itemize}
We also know that $d^2(u \otimes a) = 2(\iota \otimes c) $ because $\dot{T} S^2 = SO(3)$ and its fundamental group is $\integer_2$.

Therefore we have the following set of identities:
\begin{eqnarray*}
2 (\iota \otimes x^2) & = & d^2(u \otimes x) \\
& = & d^2 (\psi_*(u \otimes a)) \\
& = & \psi_* (d^2 ( u \otimes a)) \\
& = & \psi_* 2 (\iota \otimes c)
\end{eqnarray*}
and this implies that $\psi_*(\iota \otimes c)=\iota \otimes x^2$. Since $\iota \otimes c$ represents the class $[S^1 \times S^1]$ we can conclude that $f_*([S^1 \times S^1]) = \psi_*(\iota \otimes c) = \iota \otimes x^2 = av$.

 \end{proof}

%

\bibliographystyle{gtart}
\bibliography{Orbifoldstringbib}

\begin{thebibliography}{CJY04}

\bibitem[ACG]{AbbaspourCohenGruher}
H.~Abbaspour, R.~Cohen, and K.~Gruher.
\newblock {String topology of Poincar\'e duality groups}.
\newblock Iwase, Norio (ed.) et al., Proceedings of the conference on groups,
  homotopy and configuration spaces, University of Tokyo, Japan, July 5--11,
  2005 in honor of the 60th birthday of Fred Cohen. Coventry: Geometry \&amp;
  Topology Publications. Geometry and Topology Monographs 13, 1-10 (2008).

\bibitem[Ada78]{Adams}
J.F. Adams.
\newblock {\em Infinite loop spaces}, volume~90 of {\em Annals of Mathematics
  Studies}.
\newblock Princeton University Press, Princeton, N.J., 1978.

\bibitem[ALR]{AdemLeidaRuan}
A.~Adem, J.~Leida, and Y.~Ruan.
\newblock Orbifolds and stringy topology.
\newblock Cambridge Tracts in Mathematics 171. Cambridge: Cambridge University
  Press. xii, 149~p. 2007.

\bibitem[AR03]{AdemRuan}
A.~Adem and Y.~Ruan.
\newblock Twisted orbifold {$K$}-theory.
\newblock {\em Comm. Math. Phys.}, 237(3):533--556, 2003.

\bibitem[BH99]{BridsonHaefliger}
M.~Bridson and A.~Haefliger.
\newblock {\em Metric spaces of non-positive curvature}, volume 319 of {\em
  Grundlehren der Mathematischen Wissenschaften [Fundamental Principles of
  Mathematical Sciences]}.
\newblock Springer-Verlag, Berlin, 1999.

\bibitem[BV85]{BV}
I.~A. Batalin and G.~A. Vilkovisky.
\newblock Existence theorem for gauge algebra.
\newblock {\em J. Math. Phys.}, 26(1):172--184, 1985.

\bibitem[CG04]{CohenGodin}
R.~Cohen and V.~Godin.
\newblock A polarized view of string topology.
\newblock In {\em Topology, geometry and quantum field theory}, volume 308 of
  {\em London Math. Soc. Lecture Note Ser.}, pages 127--154. Cambridge Univ.
  Press, Cambridge, 2004.

\bibitem[CJ02]{CohenJones}
R.~Cohen and J.~Jones.
\newblock A homotopy theoretic realization of string topology.
\newblock {\em Math. Ann.}, 324(4):773--798, 2002.

\bibitem[CJY04]{CohenJonesYan}
R.~Cohen, J.~Jones, and J.~Yan.
\newblock The loop homology algebra of spheres and projective spaces.
\newblock In {\em Categorical decomposition techniques in algebraic topology
  (Isle of Skye, 2001)}, volume 215 of {\em Progr. Math.}, pages 77--92.
  Birkh\"auser, Basel, 2004.

\bibitem[CS]{ChasSullivan}
M.~Chas and D.~Sullivan.
\newblock String topology.
\newblock arXiv:math.GT/9911159.

\bibitem[CV06]{CohenVoronov}
R.~Cohen and A.~Voronov.
\newblock Notes on string topology.
\newblock In {\em String topology and cyclic homology}, Adv. Courses Math. CRM
  Barcelona, pages 1--95. Birkh\"auser, Basel, 2006.

\bibitem[Dol63]{Dold}
A.~Dold.
\newblock Partitions of unity in the theory of fibrations.
\newblock {\em Ann. of Math. (2)}, 78:223--255, 1963.

\bibitem[DW90]{DijkgraafWitten}
Robbert Dijkgraaf and Edward Witten.
\newblock Topological gauge theories and group cohomology.
\newblock {\em Comm. Math. Phys.}, 129(2):393--429, 1990.

\bibitem[Get94]{Getzler}
E.~Getzler.
\newblock Batalin-{V}ilkovisky algebras and two-dimensional topological field
  theories.
\newblock {\em Comm. Math. Phys.}, 159(2):265--285, 1994.

\bibitem[Jon87]{Jones}
J.~Jones.
\newblock Cyclic homology and equivariant homology.
\newblock {\em Invent. Math.}, 87(2):403--423, 1987.

\bibitem[LU02]{LupercioUribeLoopGroupoid}
E.~Lupercio and B.~Uribe.
\newblock Loop groupoids, gerbes, and twisted sectors on orbifolds.
\newblock In {\em Orbifolds in mathematics and physics (Madison, WI, 2001)},
  volume 310 of {\em Contemp. Math.}, pages 163--184. Amer. Math. Soc.,
  Providence, RI, 2002.

\bibitem[LU04]{LupercioUribe5}
E.~Lupercio and B.~Uribe.
\newblock Inertia orbifolds, configuration spaces and the ghost loop space.
\newblock {\em Q. J. Math.}, 55(2):185--201, 2004.

\bibitem[LU06]{LupercioUribe6}
E.~Lupercio and B.~Uribe.
\newblock Holonomy for gerbes over orbifolds.
\newblock {\em J. Geom. Phys.}, 56(9):1534--1560, 2006.

\bibitem[Moe02]{Moerdijk2002}
I.~Moerdijk.
\newblock Orbifolds as groupoids: an introduction.
\newblock In {\em Orbifolds in mathematics and physics (Madison, WI, 2001)},
  volume 310 of {\em Contemp. Math.}, pages 205--222. Amer. Math. Soc.,
  Providence, RI, 2002.

\bibitem[Seg68]{Segal1}
G.~Segal.
\newblock Classifying spaces and spectral sequences.
\newblock {\em Inst. Hautes {E}tudes Sci. Publ. Math.}, 34:105--112, 1968.

\bibitem[Vor05]{Voronov}
Alexander~A. Voronov.
\newblock Notes on universal algebra.
\newblock In {\em Graphs and patterns in mathematics and theoretical physics},
  volume~73 of {\em Proc. Sympos. Pure Math.}, pages 81--103. Amer. Math. Soc.,
  Providence, RI, 2005.

\end{thebibliography}

\end{document}